\documentclass[11pt]{article}

\usepackage{amsmath}

\usepackage{theorem}

\usepackage{amssymb}
\usepackage{latexsym}

\newcommand{\naturals}{\ensuremath{\mathbb{N}}}
\newcommand{\integers}{\ensuremath{\mathbb{Z}}}

\newcommand{\complexes}{\ensuremath{\mathbb{C}}}

\newcommand{\lin}{\mathcal{L}}
\newcommand{\bdd}{\mathcal{B}}
\newcommand{\cpct}{\mathcal{K}}
\newcommand{\tc}{\mathcal{T}}
\newcommand{\fr}{\mathcal{F}}

\newcommand{\cb}{\mathcal{CB}}

\newcommand{\seqzero}{c_{0}}
\newcommand{\seqfin}{c_{c}}

\newcommand{\isom}{\cong}

\newcommand{\cross}{\times}
\newcommand{\directsum}{\oplus}

\newcommand{\comp}{\circ}

\newcommand{\tensor}{\ensuremath{\otimes}}

\newcommand{\cmintensor}[1]{\mathop{\overset{\smile}{\text{\raisebox{0pt}[0.9ex]{\tensor}}}_{\text{#1}}}}

\newcommand{\cmaxtensor}[1]{\mathop{\overset{\frown}{\text{\raisebox{0pt}[0.7ex]{\tensor}}}_{\text{#1}}}}

\newcommand{\cBmintensor}{\cmintensor{}}
\newcommand{\cBmaxtensor}{\cmaxtensor{}}
\newcommand{\copmintensor}{\cmintensor{op}}
\newcommand{\copmaxtensor}{\cmaxtensor{op}}

\newcommand{\mintensor}[1]{\mathop{\text{\raisebox{0pt}[1ex][0ex]{$\underset{\text{\raisebox{1.1ex}[1ex][0ex]{$\smile$}}}{\otimes}$}}_{\text{#1}}}}

\newcommand{\Bmintensor}{\mintensor{}}

\newcommand{\opmintensor}{\mintensor{op}}

\newcommand{\op}{\boldsymbol}
\newcommand{\oplangle}{\langle\langle}
\newcommand{\oprangle}{\rangle\rangle}

\newcommand{\ophull}{\operatorname{\textbf{\textrm{co}}}}

\newcommand{\id}{\operatorname{id}}

\newcommand{\cmp}{\circ}
\newcommand{\clos}[1]{\overline{#1}}
\newcommand{\hull}{\operatorname{co}}
\newcommand{\spn}{\operatorname{span}}
\newcommand{\diag}{\operatorname{diag}}

\newcommand{\cone}{\operatorname{cone}}

\newtheorem{thm}{Theorem}[section]
\newtheorem{prop}[thm]{Proposition}
\newtheorem{lemma}[thm]{Lemma}
\newtheorem{cor}[thm]{Corollary}

\theoremstyle{break}
\newtheorem{defn}{Definition}[section]

\theorembodyfont{\rmfamily}
\newtheorem{example}{Example}[section]

\newenvironment{proof}[1][]%
{\par\noindent\emph{Proof#1:}\\\hspace*{\parindent}}%
{\hfill$\Box$\par\vspace{0.1in}}

\title{Matrix Compact Sets and\\
Operator Approximation Properties}

\author{Corran Webster}

\date{April 16, 1998}

\begin{document}

\maketitle

Department of Mathematics,

Texas A\&M University,

College Station,

TX 77843, USA

\vspace{2in}

* The author was supported in part by the NSF

\newpage

\section*{Abstract}
  The relationship between the operator approximation property and 
  the strong operator approximation property has deep significance 
  in the theory of operator algebras.  The original definitions of 
  Effros and Ruan, unlike the classical analogues, make no mention 
  of compact operators or compact sets.  In this paper we introduce 
  ``compact matrix sets'' which correspond to the two different 
  operator approximation properties, and show that a space has the 
  operator approximation property if and only if the ``operator 
  compact'' operators are contained in the closure of the finite 
  rank operators.  We also investigate when the two types of 
  compactness agree, and introduce a natural condition which 
  guarantees that they do.

\section{Introduction}

Since the inception of functional analysis, questions revolving around 
the approximation property have been a fruitful and important area of 
investigation.  With the development of the theory of operator spaces 
and the realization that they are ``non-commutative Banach spaces,'' 
the various versions of the approximation property have had a similar 
impact on the field, and dig to the heart of important questions in 
operator algebras.

The origin of the approximation property was the following fundamental 
result in the theory of operators on Hilbert spaces: if $X$ and $Y$ 
are Hilbert spaces, then the compact operators from $Y$ to $X$ are the 
norm closure of the finite-rank maps.  An obvious question to ask was 
whether this result extended to more general Banach spaces.  One 
direction of this---that the closure of the finite rank maps sits 
inside the compact maps---is straightforward.  The converse 
direction---that for every Banach space $Y$ every compact operator in 
$\bdd(Y,X)$ can be norm approximated by finite rank maps---is the 
\emph{approximation property} for $X$.  The initial conjecture was that every 
Banach space satisfied the approximation property.

Grothendieck~\cite{grothen:tensorprod} attacked this question, and 
although he was unable to prove or disprove the conjecture, he did 
provide a number of equivalent formulations of the approximation 
property:

\begin{thm}\label{thm:approxprop}
  Let $X$ be a Banach space.  The following are equivalent:
  \begin{enumerate}
	\item the identity map $\id: X \to X$ can be approximated 
	uniformly on compact sets by finite rank maps.,
  
	\item For all Banach spaces $Y$, the finite rank maps are dense 
	in $\bdd(Y,X)$ with the topology of uniform convergence on 
	compact sets,
  
	\item For all Banach spaces $Y$, the finite rank maps are dense 
	in $\bdd(X,Y)$ with the topology of uniform convergence on 
	compact sets,
  
	\item If $\omega \in X^{*} \cBmaxtensor X$ and $\omega(x) = 0$ 
	for all $x \in X$ then $\tau(\omega) = 0$, where $\tau(\sum 
	\chi_{i} \tensor x_{i}) = \sum \chi_{i}(x_{i})$,
  
	\item Given any Banach space $Y$, any compact map in $\bdd(Y,X)$ 
	can be approximated uniformly by finite rank maps.
  \end{enumerate}
\end{thm}

It was not until 1973 that Enflo~\cite{enflo:counterAP} provided a 
counterexample to the conjecture.  A complete discussion of the 
classical results and their implications can be found in Lindenstrauss 
and Tzafriri~\cite{lindtzaf:banachsp}.

Following Ruan's abstract characterization of concrete operator spaces 
(closed linear spaces of operators on a Hilbert space) as matrix 
normed spaces~\cite{ruan:C*subspaces}, it was realized that there were 
many parallels between the theory of Banach spaces and the theory of 
operator spaces.  The difficulty then, as now, lay in unwinding 
classical definitions and results so that they can be expressed in 
terms of concepts which have well-developed analogues in the operator 
space theory.

In their article on operator approximation properties 
\cite{effruan:approx}, Effros and Ruan observed that classically the 
following notions of convergence all agree:
\begin{enumerate}
  \item $\varphi_{\nu}$ converges to $\varphi$ uniformly on compact 
  sets in $V$,

  \item $\id \tensor \varphi_{\nu} : \seqzero \cBmintensor V \to 
  \seqzero \cBmintensor W$ converges point-norm to $\id \tensor 
  \varphi$,

  \item $\id \tensor \varphi_{\nu} : X \cBmintensor V \to X 
  \cBmintensor W$ converges point-norm to $\id \tensor \varphi$, for 
  any Banach space $X$,
  
  \item $\id \tensor \varphi_{\nu} : \ell^{\infty}(S) \cBmintensor V 
  \to \ell^{\infty}(S) \cBmintensor W$ converges point-norm to $\id 
  \tensor \varphi$ for any set $S$,
\end{enumerate}
where $V$ and $W$ are Banach spaces, $\varphi_{\nu}$ a bounded net of 
maps in $\bdd(V,W)$, $\varphi$ in $\bdd(V,W)$ and $\cBmintensor$ 
denotes the minimal Banach space tensor product.

They then noted that in the category of operator spaces there were the 
following analogues to the above:
\begin{enumerate}
  \setcounter{enumi}{1}

  \item $\id \tensor \varphi_{\nu} : \cpct \copmintensor V \to 
  \cpct \copmintensor W$ converges point-norm to $\id \tensor 
  \varphi$,
  
  \item $\id \tensor \varphi_{\nu} : X \copmintensor V \to X 
  \copmintensor V$ converges point-norm to $\id \tensor \varphi$, for 
  any operator space $X$,
  
  \item $\id \tensor \varphi_{\nu} : \bdd(H) \copmintensor V 
  \to \bdd(H) \copmintensor W$ converges point-norm to $\id \tensor 
  \varphi$ for any Hilbert space $H$,
\end{enumerate}
where $V$ and $W$ are operator spaces, $\varphi_{\nu}$ a bounded net 
of maps in $\cb(V,W)$, $\varphi$ in $\cb(V,W)$ and $\copmintensor$ 
denotes the operator space minimal tensor product.

With this observation, the following definition is natural:
\begin{defn}\label{defn:opapprox}
  $V$ has the \emph{operator approximation property} if the identity 
  map $\id : V \to V$ can be approximated by completely bounded finite 
  rank maps in the stable point norm topology.
\end{defn}

They then proved the following theorem, in direct analogy with 
Grothendieck's Theorem~\ref{thm:approxprop}:
\begin{thm}\label{thm:opapproxprop}
  Let $V$ be an operator space. The following are equivalent:
  \begin{enumerate}
    \item $V$ has the operator approximation property.
  
	\item For all operator spaces $W$, the finite rank maps in 
	$\cb(W,V)$ are dense in the stable point norm topology.
    
	\item For all operator spaces $W$, the finite rank maps in 
	$\cb(V,W)$ are dense in the stable point norm topology.
    
	\item If $\omega \in V^{*} \copmaxtensor V \subseteq \cb(V,V)$ 
	and $\omega(v) = 0$ for all $v \in V$ then $\tau(\omega) = 0$, 
	where $\tau$ is the \emph{matricial trace} $\tau(\omega) = 
	\omega(\id)$, $\id: V \to V$.
  \end{enumerate}
\end{thm}

So it appears that the topology (ii) is the ``correct'' one for the 
operator space version of the approximation property.  But what about 
the other alternatives?  It is straightforward that (iii) and (iv) are 
equivalent and imply (ii), and this was recognized by Effros and Ruan.  
That (ii) does not imply (iii) or (iv) is due to some deep results of 
Kirchberg~\cite{kirchberg:opapprox}, and it is these that lead to the 
applications in operator algebra theory.  The approximation property 
that we get when we use the topology (iii) or (iv) instead is called 
the \emph{strong operator approximation property}.

But what of the classical topology (i) that was the principal one used 
by Grothendieck?  Moreover, what of the analogue of the original 
definition of the approximation property.  Effros and Ruan saw no 
obvious operator space analogue of topology (i) or 
Theorem~\ref{thm:approxprop} (v) and they wrote~\cite{effruan:approx}:
\begin{quotation}
  It would be of considerable interest to find an analogue of 
  (i) for operator spaces. A related problem is to formulate an 
  operator space version of Grothendieck's result that a Banach 
  space $V$ has the approximation property if and only if any 
  compact operator $K: W \to V$ is a uniform limit of 
  finite rank operators.
\end{quotation}
Now that we know that the operator approximation property and the 
strong operator approximation property are in general distinct, we 
might ask if there is an analogue of topology (i) which is equivalent 
to topologies (iii) and (iv); and can we find an analogue of 
Theorem~\ref{thm:approxprop} (v) for the strong operator approximation 
property?  If we can answer these questions we should be able to shed 
light on the important question of when the operator approximation 
property and the strong operator approximation property agree.

In this paper we are able to answer all these questions completely.  
We introduce two distinct notions of a ``matrix compact set'' which 
give the results we want for the two different approximation 
properties, and develop a condition which we call 
\emph{subcoexactness} which implies the equivalence of these two types 
of matrix compactness and hence the equivalence of the two 
approximation properties.  Moreover it turns out that subcoexactness 
is a natural condition, and is related to local reflexivity.

In the next section we will answer the first of Effros and Ruan's 
questions, defining an '\emph{operator compact} matrix set and showing 
that the appropriate version of convergence on these matrix sets is 
equivalent to the stable point-norm topology of Effros and Ruan.  In 
Section 3 we review an alternative way of thinking about operator 
spaces in terms of bimodules, due to Barry Johnson, introduce a 
version of compactness in this context and prove a bimodule version of 
Grothendieck's result in the case of some special spaces.  In the 
fourth section we use these two ideas to prove the analogue of 
Grothendieck's result for the operator approximation property.  
Section 5 contains a discussion of the strong operator approximation 
case: we introduce another sort of compactness, called ``complete 
compactness'' and show that completely uniform convergence on these 
matrix sets is equivalent to the strongly stable point-norm topology; 
we discuss the implications of Kirchberg's work; and finally prove the 
analogue of Grothendieck's result for the strong operator 
approximation property.  In the final section we introduce 
\emph{coexactness} and \emph{subcoexactness}, which are properties 
dual to exactness, and show that if an operator space is subcoexact 
then our two notions of compactness agree.  We conclude by showing 
that subcoexactness is a natural condition, and is related to local 
reflexivity.

We will now introduce some notation and conventions.  Given a 
topological space $X$ we denote the continuous functions, the bounded 
continuous functions, the functions vanishing at infinity and the 
functions vanishing off compact sets by $C(X)$, $C_{b}(X)$, 
$C_{\infty}(X)$ and $C_{c}(X)$ respectively.  We denote the sequences 
vanishing at 0, the sequences which are 0 for all but finitely many 
values, the bounded sequences, absolutely summable sequences and 
square summable sequences of complex numbers as $\seqzero$ and 
$\seqfin$, $\ell^{\infty}$, $\ell^{1}$ and $\ell^{2}$ respectively.  
We will denote the standard basis in these spaces by $\{e_{k}\}$, 
where $e_{k}$ is the sequence which is zero everywhere but the $k$th 
element.  As is usual, given a Hilbert space $H$, we let $\bdd(H)$, 
$\cpct(H)$, $\tc(H)$ and $\fr(H)$ be the bounded, compact, trace-class 
and finite rank operators on $H$ respectively.  When $H$ is 
$\ell^{2}$, we will often simply write $\bdd$, $\cpct$, $\tc$ and 
$\fr$.  If $H$ is $\complexes^{n}$, we will write $M_{n}$ and 
$\tc_{n}$ for the bounded and trace-class operators respectively.  
Given a basis $\{\xi_{\nu}\}_{\nu \in \Lambda}$ for $H$, we will let 
$e_{\mu,\nu}$ be the partial isometry that takes $\spn\{\xi_{\nu}\}$ 
to $\spn\{\xi_{\mu}\}$.  If $\Gamma \subseteq \Lambda$ we denote by 
$p_{\Gamma}$ the projection from $H$ to $\spn\{\xi_{\nu} : \nu \in 
\Gamma\}$.  If $\Lambda = \naturals$, then we let $p_{n} = 
p_{\{1,\dots,n\}}$.

If $V$ is a vector space, we denote by $M_{n}(V)$ the vector space of 
$n$ by $n$ matrices with entries in $V$.  By $M_{\infty}(V)$ we mean 
the space of infinite matrices in $V$ which have only finitely many 
non-zero entries.  It will occasionally be useful to think of 
$M_{nm}(V) \isom M_{n}(M_{m}(V)) \isom M_{m}(M_{n}(V))$ as being 
indexed by tuples $(i,j)$ where $i = 1,\dots,n$ and $j = 1,\dots,m$, 
and we will denote this as $M_{n \cross m}(V)$.  We can multiply on 
the left and right by rectangular scalar matrices in the obvious way.  
Given a map $\varphi : V \to W$ between two vector spaces we define 
$\varphi_{n} : M_{n}(V) \to M_{n}(W)$ by
\[
    \varphi_{n}(v) = [\varphi(v_{i,j})].
\]
We identify $M_{n}(\lin(V,W))$ with $\lin(V,M_{n}(W))$ by mapping the 
matrix $[\varphi_{i,j}]$ to the function $v \mapsto 
[\varphi_{i,j}(v)]$.

A (non-degenerate) \emph{pairing} of two vector spaces $V$ and $W$ is 
a bilinear function
\[
    \langle \cdot,\cdot \rangle : V \cross W \to \complexes
\]
such that if $\langle v,w \rangle = 0$ for all $v \in V$, then $w = 
0$; and if $\langle v,w \rangle = 0$ for all $w \in W$, then $v = 0$.  
So each element $v \in V$ (respectively $w \in W$) determines a linear 
functional $v : W \to \complexes$ (respectively $w : V \to 
\complexes$) by
\[
    v(w) = w(v) = \langle v,w \rangle .
\]
Given such a pairing we get a \emph{matrix pairing} of $M_{n}(V)$ and 
$M_{m}(V)$ which is a map
\[
    \oplangle \cdot,\cdot \oprangle : M_{n}(V) \cross 
    M_{m}(W) \to M_{n \cross m}
\]
where the $(i,k),(j,l)$-th entry of $\oplangle v,w \oprangle$ given by 
$\langle v_{i,j},w_{k,l} \rangle$, or equivalently
\[
    \oplangle v,w \oprangle = v_{m}(w) = w_{n}(v).
\]

We will assume that the reader is familiar with the literature on 
operator spaces---an unfamiliar reader might find the following 
references useful:~\cite{effruan:opspace, paulsen:cbmapdilat, 
arveson:C*subalgebras, wittstock:HBT, ruan:C*subspaces, 
effruan:matnorm, blecher:opspdual, blechpaul:tensorprod, 
pisier:introopspfr, pisier:introopsp, pisier:exactopsp}.  We will 
denote the completely bounded maps between two operator spaces $V$ and 
$W$ by $\cb(V,W)$, and the complete and incomplete operator space 
minimal tensor products by $V \copmintensor W$ and $V \opmintensor W$, 
and the corresponding Banach space tensor products by $V \cBmintensor 
W$ and $V \Bmintensor W$.

If $A$ is a C*-algebra, with $\mathcal{V}$ and $\mathcal{W}$ 
$A$-bimodules (respectively operator $A$-bimodules), we denote the 
bounded (resp.\ completely bounded) $A$-bilinear maps from 
$\mathcal{V}$ to $\mathcal{W}$ by $\bdd_{A}(\mathcal{V},\mathcal{W})$ 
(resp.\ $\cb_{A}(\mathcal{V},\mathcal{W})$).

We define a \emph{matrix set} $\op{K} = (K_{n})$ in a vector space $V$ 
to be a family of sets $K_{n} \subseteq M_{n}(V)$ for $n \in 
\naturals$.  We call $K_{n}$ the $n$-th \emph{level} of $\op{K}$.  If 
$V$ is a vector space we will write $\op{V}$ for the matrix set 
$(M_{n}(V))$.  We say that one matrix set $\op{K}$ is a subset of 
another $\op{L}$ if $K_{n} \subset L_{n}$ for all $n \in \naturals$, 
and define intersection and union of matrix sets by taking 
intersections, or unions respectively, at each level.  If $V$ is a 
topological vector space, then we say $\op{K}$ is \emph{closed} if it 
is closed at each level, and we define \emph{open} matrix sets 
analogously.  If $\varphi : V \to W$, then given $\op{K} \subseteq 
\op{V}$ we get the matrix image of $\op{K}$,
\[
    \op{\varphi}(\op{K}) = (\varphi_{n}(K_{n})) \subseteq \op{W},
\]
and given $\op{L} \subseteq \op{W}$ we have matrix inverse image 
\[
    \op{\varphi^{-1}}(\op{L}) = (\varphi^{-1}(L_{n})) \subseteq 
\op{V}.
\]

Following~\cite{effweb:locopsp}, we say that a matrix set $\op{X}$ is 
(absolutely) \emph{matrix convex} if
\[
  \sum_{i=1}^{k} \alpha_{i} v_{i} \beta_{i} \in X_{n}
\]
whenever $v_{i} \in X_{n_{i}}$, $\alpha_{i} \in M_{n,n_{i}}$ and 
$\beta_{i} \in M_{n_{i},n}$.

Some of this paper is taken from my doctoral dissertation at UCLA, 
working under Ed Effros.  I would like to thank Ed for his guidance, 
support and help during my graduate study.  I would also like to thank 
Zhong-Jin Ruan for his many helpful comments and several key 
observations.  In particular the results concerning the strong 
operator approximation property would not have been possible without 
the definition he derived from Saar's work and the philosophical point 
that this sort of compactness corresponds classically to total 
boundedness.  My appreciation also goes out to Soren Winkler for the 
many discussions we had concerning compactness in operator spaces and 
to Vern Paulsen and Georg Schl\"{u}chtermann for their interest in my 
work and conversations we had during the Aegean conference, and to 
Marius Junge for discussions during the Operator Spaces '97 conference 
at Texas A\&M University.

\section{Operator Compactness}

The first element of our program is to find a definition of a compact 
matrix set which is suitable for our purposes.  In Webster and 
Winkler~\cite{webwink:kmt}, Webster~\cite{webster:thesis}, Le 
Merdy~\cite{lemerdy:duality}, Weaver~\cite{weaver:ncmetric} a matrix 
set $\op{X}$ is called \emph{matrix compact} if it is compact at every 
level and completely bounded in the sense that
\[
  \sup_{n} \{ \|x\| : x \in X_{n}\} < \infty
\]
Unfortunately, although it appears to be very useful for investigating 
duality theory, this definition does not give us the results we want 
with regards to the operator approximation property.  As we will see, 
it does play a role with regards to the strong operator approximation 
property.

To find the ``correct'' definition, we need to look at how we would 
prove Theorem~\ref{thm:approxprop}.  There, the critical notion is 
that if $V$ is a Banach space, compact sets can be characterized as 
being closed subsets of the closed convex hulls of sequences which 
converge to zero in $V$ (see~\cite{lindtzaf:banachsp}, 
Proposition~1.e.2).  Without loss of generality, we can take 
absolutely convex hulls instead of convex hulls.  We can think of 
sequences which converge to 0 as being elements of $\seqzero(V) = 
\seqzero \cBmintensor V$, and given an $x \in \seqzero \cBmintensor 
V$, we can define
\[
  \hull{x} = \{ v \in V : v = (\tau \tensor \id)(x), \tau \in 
  \seqfin, \|\tau\|_{1} \le 1\}.
\]
So we are saying that $K$ is compact if there is some $x \in \seqzero 
\cBmintensor V$ such that
\[
  K \subseteq \clos{\hull{x}}.
\]

Substituting $\cpct$ for $\seqzero$ at strategic places in a 
definition for Banach spaces often leads to the correct definition of 
a concept for operator spaces.  Let $V$ be an operator space and let 
$x \in \cpct(V) \isom \cpct \copmintensor V$, and define the 
\emph{absolutely matrix convex hull} of $x$ to be $\ophull{x}$ where
\[
  (\ophull{x})_{k} = \{ v \in M_{k}(V) : v = (\sigma \tensor \id)(x), 
  \sigma \in M_{k}(M_{\infty}), \|\sigma\|_{\tc} \le 1 \}.
\]
and we are thinking of $M_{\infty}$ as sitting inside $\tc$.  This is 
an absolutely matrix convex set, since if $v \in (\ophull{x})_{k}$, $w 
\in (\ophull{x})_{k}$, then $v = (\sigma \tensor \id)(x)$, $w = (\tau 
\tensor \id)(x)$ and
\[
  v \directsum w = (\sigma \directsum \tau \tensor \id)(x)
\]
where $\sigma \directsum \tau \in M_{k+l}(M_{\infty})$ and the unit 
ball of $\tc$ is absolutely matrix convex.  Similarly
\[
  \alpha v \beta = (\alpha \sigma \beta \tensor \id)(x)
\]
and $\alpha \sigma \beta \in M_{l}(M_{\infty})$ and again, the unit 
ball of $\tc$ is absolutely matrix convex.

We also note that if $\sigma \in M_{k}(\tc)$, then $(\sigma \tensor 
\id)(x) \in \clos{\ophull(x)}$, since if we let $\sigma_{n}(v) = 
\sigma(p_{n}vp_{n})$, then $p_{n}vp_{n} \to v$ as $n \to \infty$ and 
so $(\sigma_{n} \tensor \id)(x) \to (\sigma \tensor \id)(x)$ as $n \to 
\infty$.

\begin{defn}
  If $\op{K}$ is a matrix subset of an operator space $V$, then we say 
  $\op{K}$ is \emph{operator compact} if $\op{K}$ is closed and there 
  is some $x \in \cpct(V)$ such that $\op{K} \subseteq 
  \clos{\ophull{x}}$.
\end{defn}

\begin{example}\label{eg:tcball}
Consider the matrix unit ball of $\tc_{n}$.  Let $\tau \in 
M_{n}(\tc_{n}) \isom \cb(M_{n},M_{n})$ be the identity map, or 
equivalently
\[
  \tau = \begin{bmatrix} \tau_{1,1} & \cdots & \tau_{1,n} \\
  \vdots & \ddots & \vdots \\ \tau_{n,1} & \cdots & \tau_{n,n} 
\end{bmatrix}
  = \sum_{i,j = 1}^{n} e_{i,j} \tensor \tau_{i,j}
\]
where $\tau_{i,j}(a) = a_{i,j}$.  Then given any $\sigma$ in the 
matrix unit ball of $\tc_{n}$, we have that
\[
  \sigma(a) = \sigma(\tau(a)) = (\sigma \tensor \id)(\sum_{i,j = 
  1}^{n} e_{i,j} \tensor \tau_{i,j})(a)
\]
and so $\sigma \in \ophull(\tau)$.  Hence this matrix unit ball is 
operator compact.
\end{example}

Our claim is that operator compactness is the correct definition to 
solve the questions posed by Effros and Ruan.

\begin{defn}
  Let $V$, $W$ be operator spaces, $\varphi_{\nu}$, $\varphi 
  \in \cb(V,W)$.  We say that the net $\varphi_{\nu}$ converges 
  to $\varphi$ \emph{completely uniformly on operator compact 
  sets} if for all operator compact sets $\op{X} \subset 
  V$ and $\varepsilon > 0$ there is an $N$ such that 
  \[
	\sup \{\| (\varphi_{\nu})_{n}(x) - \varphi_{n}(x) \|_{n} : 
	x \in X_{n},\ n \in \naturals \} < \varepsilon,
  \]
  for all $\nu \ge N$.
\end{defn}

Our aim is to show that this topology is the same as the topology 
(ii), which we call the \emph{stable point-norm topology}.

\begin{prop}\label{prop:convergence}
  Let $V$, $W$ be operator spaces, $\varphi_{\nu}$, $\varphi \in 
  \cb(V,W)$.  Then the following are equivalent:
  \begin{enumerate}
	\item $\varphi_{\nu} \to \varphi$ in the stable point norm 
	topology on $\cb(V,W)$.
      
	\item $\varphi_{\nu} \to \varphi$ completely uniformly on operator 
	compact sets of $V$.
  \end{enumerate}
\end{prop}

\begin{proof}
  (i $\implies$ ii): Let $\op{X}$ be any operator compact set, and 
$x 
  \in \cpct(V)$ be so that $\op{X} \subseteq \clos{\ophull{x}}$.  
We 
  know that given any $\varepsilon$ then for $\nu$ sufficiently 
large we 
  have
  \[
    \|(\varphi_{\nu} \tensor \id)(x) - (\varphi \tensor \id)(x)\| 
\le 
    \frac{\varepsilon}{2\|x\|}
  \]
  and so we have
  \begin{align*}
	&\sup\{\|\op{\varphi_{\nu}}(v) - \op{\varphi}(v)\| : v \in 
	\ophull{x}\} \\
	&= \sup\{\| (\id \tensor \varphi_{\nu})(\sigma 
	\tensor \id)(x) - (\id \tensor \varphi_{\nu})(\sigma \tensor 
	\id)(x)\| : \sigma \in M_{k}(M_{\infty}),\ \|\sigma\|_{\tc} \le 1 
\}\\
	&= \sup\{\|(\sigma \tensor \id)(\id \tensor \varphi_{\nu})(x) - 
	(\sigma \tensor \id)(\id \tensor \varphi)(x)\| : \sigma \in 
	M_{k}(M_{\infty}),\ \|\sigma\|_{\tc} \le 1 \}\\
    &\le \|(\varphi_{\nu} \tensor \id)(x) - (\varphi \tensor 
\id)(x)\| \\
    &\le \varepsilon/2.
  \end{align*}
  
  Taking closures we then get that for $\nu$ sufficiently large we 
  have
  \[
	\sup_{v \in \op{X}}\|\op{\varphi_{\nu}}(v) - 
	   \op{\varphi}(v)\|_{\infty} < \varepsilon.
  \]

  (ii $\implies$ i): Given $x \in \cpct(V)$ we notice that $\pi_{n}(x) 
  \in \ophull{x}$ for all $n$, that $\ophull{x}$ is operator compact, 
  and so if we choose $\nu$ sufficiently large we have
  \begin{equation*}\begin{split}
	\|(\varphi_{\nu} \tensor \id)(x) - (\varphi \tensor 
	    \id)(x)\|_{\infty} &= \sup_{n} \|p_{n}(\varphi_{\nu} \tensor 
	    \id)(x) - (\varphi \tensor \id)(x)p_{n}\| \\
	&= \sup_{n} \|(\varphi_{\nu} \tensor \id)(x_{n}) - (\varphi 
	    \tensor \id)(x_{n}) \\	
    &< \varepsilon.
  \end{split}\end{equation*}
\end{proof}

So we have found a simple answer to the first question posed by Effros 
and Ruan.  To prove the analogue of Grothendieck's result, however, we 
will need to look much more deeply.  We define an \emph{operator 
compact map} $\varphi$ between two operator spaces $V$ and $W$ as 
being one for which the image of the matrix unit ball is compact.

Completely bounded finite rank maps are operator compact, for as we 
will see in Corollary~\ref{cor:finitematrixball} the unit balls of 
finite dimensional operator spaces are operator compact.  Although it 
is clear that the operator compact maps are an ideal under composition 
in the completely bounded maps, it is not clear whether the operator 
compact maps are closed in the completely uniform topology.  The 
difficulty is that there is little control over the $x$ which we take 
the matrix convex hull of---and indeed, as we will see in the final 
section, in some cases there is little prospect of gaining any 
control.

\section{C*-Operator Spaces}

The missing technology that we will need is the concept of a 
C*-operator space, introduced by Barry 
Johnson~\cite{johnson:kopspace}.  If the C*-algebra is $\cpct$, then 
these spaces form a category which is equivalent to that of operator 
spaces, but has the advantage over operator spaces that often one can 
adapt Banach space methods to work for C*-operator spaces, where it 
would not necessarily even be possible to formulate a strategy for 
operator spaces.  In particular, this will be crucial in our analysis 
of the operator approximation property.

If $A$ is a C*-algebra then Johnson defines an \emph{$A$-operator 
space} $\mathcal{V}$ to be an essential $A$-bimodule with a norm which 
is absolutely $A$-convex, i.e.  if $v_{1}$, $v_{2}$ lie in the unit 
ball of $\mathcal{V}$ then so does
\begin{equation}
  v = \alpha_{1} v_{1} \beta_{1} + \alpha_{2} v_{2} \beta_{2}
\end{equation}
where $\| \alpha_{1}\alpha_{1}^{*} + \alpha_{2}\alpha_{2}^{*} \| \le 
1$ and $\| \beta_{1}^{*}\beta_{1} + \beta_{1}^{*}\beta_{1}\| \le 1$.

The natural morphisms are the continuous bi-$A$-linear maps, i.e.  
continuous maps $\varphi : \mathcal{V} \to \mathcal{W}$ such that 
$\varphi(\alpha v \beta) = \alpha \varphi(v) \beta$ for all $v \in V$, 
$\alpha$, $\beta \in A$. We will denote the space of all such maps 
by $\bdd_{A}(\mathcal{V},\mathcal{W})$.

A recent result of Magajna~\cite[Theorem 2.1]{magajna:minopmod} tells 
us that an $A$-bimodule has an operator bimodule structure---that is 
the bimodule has an operator space structure and the bimodule action 
is completely contractive---if and only if
\[
  \| \alpha_{1} v_{1} \beta_{1} + \alpha_{2} v_{2} \beta_{2} \| \le
  \| \alpha_{1}\alpha_{1}^{*} + \alpha_{2}\alpha_{2}^{*} \|^{1/2} 
  \max \{\|v_{1}\|,\|v_{2}\|\} \| \beta_{1}^{*}\beta_{1} + 
  \beta_{1}^{*}\beta_{1}\|^{1/2}.
\]
In other words, every C*-operator space has an operator bimodule 
structure (indeed, it potentially has many). In particular Magajna 
considers the \emph{minimal} operator $A$-bimodule whose norms are 
given by
\[
  \| v \|_{n} = \sup \{ \| \alpha v \beta \| : \alpha \in M_{1,n}(A), 
  \beta \in M_{n,1}(A), \|\alpha\| \le 1, \| \beta \| \le 1 \}.
\]
Given a C*-operator space $\mathcal{V}$, we denote the corresponding 
minimal operator bimodule by $\min \mathcal{V}$.  Magajna does not 
look at morphisms explicitly, however the following lemma is immediate 
from his definition, and is natural given the theory of operator 
spaces.

\begin{lemma}
  If $\mathcal{V}$ and $\mathcal{W}$ are $A$-operator spaces, then 
  $\varphi \in \bdd_{A}(\mathcal{V},\mathcal{W})$ if and only if 
  $\varphi \in \cb_{A}(\min \mathcal{V},\min \mathcal{W})$.  Indeed 
  the completely bounded norm is equal to the bounded norm.
\end{lemma}

\begin{proof}
  We simply note that
  \begin{align*}
	\| \varphi_{n}(v) \|_{n} &= \sup \{ \| \alpha \varphi_{n}(v) \beta 
	\| : \\
    &\qquad \alpha \in M_{1,n}(A), \beta \in M_{n,1}(A), \|\alpha\| \le 
	1, \| \beta \| \le 1 \} \\
	&= \sup \{ \| \varphi(\alpha v \beta) \| : \\
	&\qquad \alpha \in M_{1,n}(A), \beta \in M_{n,1}(A), \|\alpha\| 
	\le 1, \| \beta \| \le 1 \} \\
    &\le \| \varphi \| \sup \{ \| \alpha v \beta \| : \\
    &\qquad \alpha \in M_{1,n}(A), \beta \in M_{n,1}(A), \|\alpha\| \le 1, 
	\| \beta \| \le 1 \} \\
    &= \| \varphi \| \|v\|_{n},
  \end{align*}
  so $\| \varphi \| \le \| \varphi \|_{\infty} \le \| \varphi \|$. 
\end{proof}

If $A$ is an injective C*-algebra, then this immediately gives us a 
Hahn-Banach theorem for the category of $A$-operator spaces.

\begin{prop}
  Let $A$ be an injective C*-algebra.  Then if $\mathcal{V}$, 
  $\mathcal{W}$ are $A$-operator spaces such that $\mathcal{V}$ is a 
  bi-$A$-invariant subspace of $\mathcal{W}$, then given any 
  continuous bi-$A$-linear functional
  \[
    \varphi : \mathcal{V} \to A
  \]
  then there is a continuous bi-$A$-linear functional
  \[
    \bar{\varphi} : \mathcal{W} \to A
  \]
  such that $\bar{\varphi}|_{\mathcal{V}} = \varphi$ and $\| 
  \bar{\varphi} \| = \| \varphi \|$.
\end{prop}

\begin{proof}
  We can consider $\varphi \in \cb_{A}(\min \mathcal{V}, \min A)$, and 
  use Wittstock's Hahn-Banach theorem~\cite[Theorem 3.1]{wittstock:HBT} 
  for operator bimodules to find a completely bounded extension
  \[
    \bar{\varphi} : \mathcal{W} \to A.
  \]
  But by the previous lemma $\bar{\varphi}$ is also a bounded 
  extension of $\bar{\varphi}$ with $\| \bar{\varphi} \| = \| 
  \bar{\varphi} \|_{\infty} = \| \varphi \|_{\infty} = \| \varphi \|$.
\end{proof}

In particular, we shall use this result with $A = \bdd(H)$.  However, 
we will also need a Hahn-Banach theorem in the case where $A = \cpct$, 
and $\cpct$ is not injective.  Fortunately we have another approach 
which delivers us such a theorem.  In this special case we can 
rephrase our convexity axiom in much more familiar language.  A set 
$X$ in a $\cpct$-operator space is $\cpct$-convex if and only if it 
satisfies
\begin{alignat}{2}
  v + w &\in X && \qquad\text{for all orthogonal $v$, $w \in 
  X$,}\tag{A$\cpct$C1}\label{AKC1} \\
  \alpha v \beta &\in X &&\qquad \text{for all $v \in X$, $\alpha$, 
  $\beta \in \cpct$, $\| \alpha \|$, $\|\beta\| \le 1$.}
  \tag{A$\cpct$C2}\label{AKC2}
\end{alignat}
We say that $v_{1},\dots,v_{n} \in \mathcal{V}$ are \emph{orthogonal} 
if there exist orthogonal projections $e_{1},\dots,e_{n} \in \cpct$ 
such that $e_{i}v_{i}e_{i} = v_{i}$.  A $\cpct$-norm is a norm 
$\|\cdot\|$ which satisfies
\begin{alignat}{2}
  \| v + w \| &= \max\{\|v\|,\|w\|\} &&\qquad \text{for all orthogonal 
  $v$, $w \in X$,}\tag{A$\cpct$N1}\label{AKN1} \\
  \|\alpha v \beta\| &\le \|\alpha\|\|v\|\|\beta\| &&\qquad \text{for 
  all $v \in X$, $\alpha$, $\beta \in \cpct$.} 
  \tag{A$\cpct$N2}\label{AKN2}
\end{alignat}
Clearly the unit balls of $\cpct$-norms are $\cpct$-convex.

Given the parallels between these axioms and Ruan's axioms for the 
matrix norms of operator spaces, it is not surprising that $\cpct 
\copmintensor V$ is a $\cpct$-operator space where $\cpct$ operates 
via
\[
  \alpha(a \tensor v)\beta = \alpha a \beta \tensor v.
\]
Moreover the natural map $\varphi \mapsto \id \tensor \varphi = 
\varphi_{\infty}$ is an isometric isomorphism between completely 
bounded maps and continuous bi-$\cpct$-linear maps on the 
corresponding $\cpct$-operator spaces.

Hence $\cpct : V \mapsto \cpct(V)$ takes operator spaces to 
$\cpct$-operator spaces, and $\cpct : \varphi \mapsto 
\varphi_{\infty}$ takes completely bounded maps to continuous 
bi-$\cpct$-linear maps.  In other words $\cpct$ is a functor.

We define
\[
  \fr(\mathcal{V}) = \{ e_{1,1} v e_{1,1} : v \in \mathcal{V} \}.
\]
Then $V = \fr(\mathcal{V})$ is an operator space when we identify 
$M_{n}(V)$ with
\[
  \{ p_{n} v p_{n} : v \in \mathcal{V} \}
\]
via the map
\[
  \tau_{n} : [v_{i,j}] \mapsto \sum_{i,j=1}^{n} e_{i,1} v_{i,j} 
e_{1,j},
\]
and give it the norm inherited from $\mathcal{V}$.  That this is in 
fact an operator space norm follows immediately from (\ref{AKN1}) and 
(\ref{AKN2}).  If
\[
  \varphi : \mathcal{V} \to \mathcal{W}
\]
is a continuous bi-$\cpct$-linear then $\varphi|_{V}$ is a completely 
bounded map from $V$ to $W = \fr(\mathcal{W})$, since
\[
  (\varphi|_{V})_{n}([v_{i,j}]) = 
  \varphi|_{p_{n}\mathcal{V}p_{n}}(\tau_{n}([v_{i,j}]))
\]
and again, this is an isometric isomorphism of the appropriate 
mapping spaces.  Hence $\fr$ is a functor from $\cpct$-operator 
spaces to operator spaces.

The next proposition, to the author's knowledge first explicitly 
stated by Johnson, but not explicitly proved, is now part of the 
folklore of the subject.

\begin{prop}
  The functors $\cpct$ and $\fr$ implement an equivalence of 
  categories between the category of operator spaces with completely 
  bounded maps and the category of $\cpct$-operator spaces with 
  continuous bi-$\cpct$-linear maps.
\end{prop}

The proof is essentially a matter of resolving semantic differences, 
and is not of particular interest.  For this reason we relegate it to 
an Appendix.  Note that there are many possible ways of implementing 
this equivalence.  The point of such a formal result, is that we can 
quickly transfer results from the theory of operator spaces to 
$\cpct$-operator spaces and vice-versa.  In particular, it allow us to 
quickly prove a Hahn-Banach theorem for $\cpct$-operator spaces.

\begin{prop}
  If $\mathcal{V}$, $\mathcal{W}$ are $\cpct$-operator spaces such 
  that $\mathcal{V}$ is a bi-$\cpct$-invariant subspace of 
  $\mathcal{W}$, then given any continuous bi-$\cpct$-linear 
  functional
  \[
    \varphi : \mathcal{V} \to \cpct
  \]
  then there is a continuous bi-$\cpct$-linear functional
  \[
    \bar{\varphi} : \mathcal{W} \to \cpct
  \]
  such that $\bar{\varphi}|_{\mathcal{V}} = \varphi$ and $\| 
  \bar{\varphi} \| = \| \varphi \|$.
\end{prop}

\begin{proof}
  We have shown that we can find $\psi = \fr(\varphi)$ which is a 
  completely bounded linear functional on $V = \fr(\mathcal{V})$.  The 
  bi-$\cpct$-invariance of $\mathcal{V}$ in $\mathcal{W}$ implies that 
  $V$ is a subspace of $W = \fr(\mathcal{W})$, and so we get by the 
  Hahn-Banach theorem for operator spaces that there is a completely 
  bounded linear functional $\bar{\psi}$ on $W$ extending $\psi$.  We 
  then push this back to the original category to get a 
  bi-$\cpct$-linear functional $\bar{\varphi}$ which extends 
  $\varphi$.  It remains only to note that the norms are preserved by 
  the functors.
\end{proof}

We are now in a position to discuss the appropriate notion of 
compactness in these categories.

\begin{defn}
  Let $\mathcal{V}$ be an $A$-operator space. We say that $X \subset 
  \mathcal{V}$ is \emph{$A$-compact} if $X$ is closed and
  \[
    X \subseteq \clos{\hull_{A}} \{ x_{i} \}
  \]
  where $x_{i} \to 0$ and $\clos{\hull_{A}}$ indicates the closed 
  $A$-convex hull.

  We say that a bounded bi-$A$-linear map $\varphi : \mathcal{V} \to 
  \mathcal{W}$ is an \emph{$A$-compact} operator if the image of the 
  unit ball of $\mathcal{V}$ is $A$-compact in $\mathcal{W}$.
\end{defn}

We now have analogues of Proposition~\ref{prop:convergence}:

\begin{prop}\label{prop:convergence2}
  Let $\mathcal{V}$, $\mathcal{W}$ be $A$-operator spaces, 
  $\varphi_{\nu}$, $\varphi \in \bdd_{A}(\mathcal{V},\mathcal{W})$.  
  Then the following are equivalent:
  \begin{enumerate}      
	\item $\id \tensor \varphi_{\nu} \to \id \tensor \varphi$ 
	point-norm in 
	$\bdd_{A}(\seqzero(\mathcal{V}),\seqzero(\mathcal{W}))$.
      
	\item $\varphi_{\nu} \to \varphi$ uniformly on the compact sets of 
	$\mathcal{V}$.
  \end{enumerate}
\end{prop}

\begin{proof}
  (i $\implies$ ii): If $X$ is any $A$-compact set, $X \subseteq 
  \clos{\hull_{A}{\{x_{i}\}}}$, then
  \begin{align*}
	\sup &\{ \| \varphi_{\nu}(v) - \varphi(v) \| : v \in 
	\hull_{A}{\{x_{i}\}}\} \\
	&\le \sup_{n \in \naturals} \{ \| \varphi_{\nu}(\sum_{i=1}^{n} 
	\alpha_{i}x_{i}\beta_{i}) - \varphi(\sum_{i=1}^{n} 
	\alpha_{i}x_{i}\beta_{i}) \| : \|\sum_{i=1}^{n} 
	\alpha_{i}^{*}\alpha_{i}\|,\ \|\sum_{i=1}^{n} 
	\beta_{i}\beta_{i}^{*}\| \le 1\} \\
	&\le \sup _{n \in \naturals} \{ \| \sum_{i=1}^{n} \alpha_{i} 
	(\varphi_{\nu}(x_{i}) - \varphi(x_{i})) \beta_{i} \| : 
	\|\sum_{i=1}^{n} \alpha_{i}^{*}\alpha_{i}\|,\ \|\sum_{i=1}^{n} 
	\beta_{i}\beta_{i}^{*}\| \le 1\} \\
	&\le \sup \{ \| \varphi_{\nu}(x_{i}) - \varphi(x_{i}) \| \}
  \end{align*}
  since balls are $A$-convex.  Taking closures then gives the 
  result.

  (ii $\implies$ i): Sequences converging to 0 are $A$-compact, and so 
  the result is immediate.
\end{proof}

We would like to define an analogue of the operator approximation 
property in these categories, but there is an obstacle revolving 
around what is the appropriate analogue of a finite rank map.  Our 
interest in later sections is going to concentrate on $A$-operator 
spaces of the form $A \copmintensor V$ for some operator space $V$.  
As a result we say that a bi-$A$-linear map
\[
  \varphi : \mathcal{W} \to \mathcal{V}
\]
has \emph{finite $A$-rank} if the range is bi-$A$-linearly 
isometrically isomorphic to $A \copmintensor E$, where $E$ is a finite 
dimensional operator space. In particular if $\mathcal{V} = A 
\copmintensor V$, then $\varphi$ may be written as
\[
  \varphi = \sum_{i=1}^{k} \varphi_{i} \tensor v_{i}
\]
where $v_{i} \in V$ and $\varphi_{i} \in \bdd_{A}(\mathcal{W},A)$.

\begin{defn}
  We say that an $A$-operator space $\mathcal{V}$ has the 
  $A$-approximation property if the identity map can be approximated 
  uniformly on $A$-compact sets by finite $A$-rank maps in 
  $\bdd_{A}(\mathcal{V},\mathcal{V})$.
\end{defn} 

With this definition and \ref{prop:convergence2} we have the 
analogous version of Effros and Ruan's result.

\begin{thm}
  The following conditions are equivalent to a $A$-operator space 
  $\mathcal{V}$ having the $A$-operator approximation property:
  \begin{enumerate}
	\item For all $A$-operator spaces $\mathcal{W}$, the finite 
	$A$-rank maps are dense in $\bdd_{A}(\mathcal{W},\mathcal{V})$ 
	with topology of uniform convergence on $A$-compact sets.
    
	\item For all $A$-operator spaces $\mathcal{W}$, the finite 
	$A$-rank maps are dense in $\bdd_{A}(\mathcal{V},\mathcal{W})$ 
	with topology of uniform convergence on $A$-compact sets.
  \end{enumerate}
\end{thm}

\begin{proof}
  Clearly either of these implies the $A$-operator approximation 
  property as a special case.

  On the other hand if finite $A$-rank maps $\varphi_{\nu}$ converge 
  to $\id$, then for any $\varphi \in 
  \bdd_{A}(\mathcal{W},\mathcal{V})$ (resp.\ 
  $\bdd_{A}(\mathcal{V},\mathcal{W})$) the maps $\varphi_{\nu} \circ 
  \varphi$ (resp.\ $\varphi \circ \varphi_{\nu}$) are finite $A$-rank 
  maps. But
  \[
    \varphi_{\nu} \circ \varphi \to \varphi
  \]
  and
  \[
    \varphi \circ \varphi_{\nu} \to \varphi
  \]
  uniformly on $A$-compact sets
\end{proof}

We now have all the pieces to prove a version of Theorem 1.1 
in the context of certain C*-operator spaces.  The spaces we work with 
are determined by the restriction on those categories of C*-operator 
spaces for which we have a Hahn-Banach theorem available, and 
consideration of what we need the results for in the sequel.  More 
general results may be possible.

\begin{thm}\label{thm:C*-approx}
  Let $A$ be either an injective C*-algebra or $\cpct$.  Let 
  $\mathcal{V} = A \copmintensor V$ for some operator space $V$.  Then 
  $\mathcal{V}$ has the $A$-approximation property if and only if for 
  any $A$-operator space $\mathcal{W}$ and any $A$-compact map 
  $\varphi \in \bdd_{A}(\mathcal{W},\mathcal{V})$, the map $\varphi$ 
  can be approximated uniformly by $A$-finite rank maps.
\end{thm}

The proof is really just an adaptation of the classical proof to a 
bi-module setting.

\begin{proof}
  First assume that $\mathcal{V}$ has the $A$-approximation property.  
  Since $\varphi$ is $A$-compact, there is a sequence $x = (x_{i}) \in 
  \seqzero \cBmintensor \mathcal{V}$ such that $\varphi(B) \subseteq 
  \clos{\hull_{A}}\{ x_{i} \}$, where $B$ is the unit ball of 
  $\mathcal{W}$.  Now we know that for any $\varepsilon > 0$ we can 
  find a finite rank map $\psi \in \bdd_{A}(\mathcal{V},\mathcal{V})$ 
  such that
  \[
    \|\psi(v) - v\| < \varepsilon
  \]
  for all $v \in \clos{\hull_{A}}\{x\}$ and so
  \[
    \| \psi \cmp \varphi - \varphi \| < \varepsilon.
  \]

  Conversely, let $X$ be any $A$-compact set in $\mathcal{V}$, and 
  without loss of generality, we may assume that
  \[
	X = \clos{\hull_{A}} \{ x_{i} \}
  \]
  for some sequence $\{ x_{i} \} \in \seqzero(\mathcal{V})$, and that 
  $x_{i} \ne 0$ for all $i$.  Now let
  \[
	U = \clos{\hull_{A}} \left\{ \frac{x_{i}}{\|x_{i}\|^{1/2}} 
	\right\},
  \]
  and so the identity map 
  \[
    \id : {}_{U}\mathcal{V} \to V,
  \]
  where ${}_{U}\mathcal{V}$ is the vector space $\mathcal{V}$ with the 
  $A$-convex norm determined by $U$, is $A$-compact.  Therefore we can 
  find finite $A$-rank maps
  \[
    \varphi_{k} : {}_{U}\mathcal{V} \to V
  \]
  which approximate $\id$ in norm; in particular, we may let each map 
  be of the form
  \[
	\varphi_{k}(v) = \sum_{i = 1}^{n_{k}} \psi_{k,i}(v) \tensor v_{k,i}
  \]
  where $\psi_{k,i} \in \bdd_{A}({}_{U}\mathcal{V},A)$.  So all we 
  need do is show we can approximate elements of 
  $\bdd_{A}({}_{U}\mathcal{V},A)$ by elements of 
  $\bdd_{A}(\mathcal{V},A)$ uniformly on $X$, i.e. given $\delta > 0$, 
  then if for all $k$, $i$ we can find $\psi_{k,i}' \in 
  \bdd_{A}(\mathcal{V},A)$ such that
  \[
	\| \psi_{k,i}'(v) - \psi_{k,i}(v) \| < \delta
  \]
  for all $v \in X$, then if we let
  \[
	\varphi_{k}'(v) = \sum_{i = 1}^{n_{k}} \psi_{k,i}'(v) \tensor v_{k,i}
  \]
  we have
  \[
	\| \varphi_{k}'(v) - v \| \le \| \varphi_{k}'(v) - 
	\varphi_{k}(v) \| + \| \varphi_{k}(v) - v \| < 2\delta
  \]
  for all $k$ sufficiently large, whence the result.

  So given any $\psi \in \bdd_{A}({}_{U}\mathcal{V},A)$ 
  we may assume $\| \psi \|_{U} = 1$ without loss of 
  generality.  We note that $\frac{x_{i}}{\|x_{i}\|^{1/2}} 
  \in U$ for all $i$, and so we have that
  \[
    \left\| \frac{x_{i}}{\|x_{i}\|^{1/2}} \right\|_{U} \le 1
  \]
  for all $i$.  But this is equivalent to saying that
  \[
	\| x_{i} \|_{U} \le \| x_{i} \|^{1/2}
  \]
  for all $i$, and so $\| x_{i} \|_{U} \to 0$ as $i \to \infty$.  
  Now choose an $N$ so that $\| x_{i} \| < \delta^{2}/2$ for all $i 
  > N$.  Then let
  \[
	\mathcal{V}_{\delta} = \spn_{A}{\{ x_{i} \}_{i=1}^{N}},
  \]
  and let $\psi_{\delta}$ be the restriction of $\psi$ to 
  $\mathcal{V}_{\delta}$.  $\psi_{\delta}$ is in 
  $\bdd_{\cpct}(\mathcal{V}_{\delta},\cpct)$ since 
  $\mathcal{V}_{\delta}$ is a finite bi-$A$-linear span.  On the other 
  hand
  \[
    \|\psi_{\delta}\| \le \|\psi_{\delta}\|_{U} = 1,
  \]
  since $U$ sits inside the unit ball of $\mathcal{V}$.  But by the 
  Hahn-Banach theorem for $A$-operator spaces, we can find a 
  $\psi' \in \bdd_{A}(\mathcal{V},A)$ which agrees with 
  $\psi_{\delta}$ on $\mathcal{V}_{\delta}$ and
  \[
	\| \psi' \| = \| \psi_{\delta} \| \le 1.
  \]

  Now for any $x_{i}$ we have either that $\psi'(x_{i}) = \psi(x_{i})$ 
  (if $i \le N$), or
  \begin{eqnarray*}
    \| \psi(x_{i}) - \psi'(x_{i}) \|
        & \le & \| \psi(x_{i}) \| + \| \psi'(x_{i}) \| \\
        & \le & \| \psi \| \| x_{i} \|_{U} + \| \psi'\| \| x_{i}\| \\
        & \le & \delta/2 + \delta^{2}/2 \le \delta\ \text{(for $\delta$ small).}
  \end{eqnarray*}
  and it is easy to see that this implies $\| \psi(x) - \psi'(x) \| < 
  \delta$ for all $x \in X$.
\end{proof}

\section{The Operator Approximation Property}

The point of the previous section was to introduce notions of 
compactness which allows us to use more classical techniques. What we 
must now do is to relate those notions with the ideas from Section 2 
so that we can get the result we need. Our first step is to relate 
Propositions~\ref{prop:convergence} and~\ref{prop:convergence2} by 
showing that across the equivalence of categories that the convergence 
is the same.

\begin{prop}\label{prop:convergence3}
  Let $V$, $W$ be operator spaces, $\varphi_{\nu}$, $\varphi \in 
  \cb(V,W)$.  Then the following are equivalent:
  \begin{enumerate}
	\item $\varphi_{\nu} \to \varphi$ in the stable point norm 
	topology on $\cb(V,W)$.
      
	\item $\varphi_{\nu} \to \varphi$ completely uniformly on operator 
	compact sets of $V$.
      
	\item $\id \tensor (\varphi_{\nu})_{\infty} \to \id \tensor 
	\varphi_{\infty}$ point-norm in 
	$\bdd_{\cpct}(\seqzero(\cpct(V)),\seqzero(\cpct(V)))$.
      
	\item $(\varphi_{\nu})_{\infty} \to \varphi_{\infty}$ uniformly on 
	the compact sets of $\cpct(V)$.
      
	\item $(\varphi_{\nu})_{\infty} \to \varphi_{\infty}$ uniformly on 
	the $\cpct$-compact sets of $\cpct(V)$.
  \end{enumerate}
\end{prop}

However to prove this, we will need the following lemma about tensor 
products of commutative C*-algebras with operator spaces.

\begin{lemma}\label{lemma:tensor}
  If $V$ is an operator space, and $X$ is a locally compact Hausdorff 
  space, then $C_{0}(X) \cBmintensor V \isom C_{0}(X) \copmintensor V$ 
  as Banach spaces.
\end{lemma}

\begin{proof}[ (Lemma~\ref{lemma:tensor})]
  Recall that there exists a C*-algebra $A$ such that $V$ embeds 
  completely isometrically in $A$, and that $C_{0}(X) \cBmintensor A 
  \isom C_{0}(X) \copmintensor A$ as Banach spaces.  However, the 
  minimal tensor products respect inclusion, and so $C_{0}(X) 
  \Bmintensor V$ and $C_{0}(X) \opmintensor V$ are isometrically 
  isomorphic as normed vector spaces, and so their completions 
  agree.
\end{proof}

\begin{proof}[ (Proposition~\ref{prop:convergence2})]
  (i $\iff$ ii): Proposition~\ref{prop:convergence}
  
  (i $\implies$ iii): We observe that
  \[
	\cpct(V) \isom \cpct \copmintensor V \isom \cpct \copmintensor 
	\cpct \copmintensor V \supseteq c_{0} \copmintensor \cpct 
	\copmintensor V \isom c_{0} \cBmintensor (\cpct \copmintensor V)
  \]
  as Banach spaces by Lemma~\ref{lemma:tensor}.  So if 
  $\varphi_{\nu} \to \varphi$ in the stable point norm topology, 
  then
  \[
	\id \tensor \id \tensor \varphi_{\nu} \to \id \tensor \id \tensor 
	\varphi
  \]
  in the point-norm topology on $c_{0} \cBmintensor \cpct(V)$.

  (iii $\iff$ iv): This is the classical result, since $\cpct(V)$, 
  $\cpct(W)$ are Banach spaces.

  (iv $\implies$ i): Follows since a point in $\cpct(V)$ is compact.

  (iii $\iff$ v): Proposition~\ref{prop:convergence2} with $A = \cpct$.
\end{proof}

\begin{cor}\label{cor:approxpropequiv}
  An operator space $V$ has the operator approximation property if and 
  only if $\cpct(V)$ has the $\cpct$-approximation property.
\end{cor}

We still need a little more, however, since we need to relate operator 
compact operators to $\cpct$-compact operators so that we can use 
Theorem~\ref{thm:C*-approx}. In fact we can show that $\cpct$-compact 
convex sets correspond to sets which are operator compact.

To see this, we assume that
\[
  X \subseteq \clos{\hull_{\cpct}(\{x_{1},x_{2},\dots\})}.
\]
Then we let $\lambda$ be an isomorphism of $\cpct \cBmintensor 
\cpct$ with $\cpct$, where, for simplicity's sake, we will assume 
$\lambda$ is induced by a bijection $\mu: \naturals \cross 
\naturals \to \naturals$ via
\[
  \lambda(e_{(i,j),(k,l)}) = e_{\mu(i,j),\mu(k,l)}
\]
and let $x = \lambda (\diag(x_{1},x_{2},\dots))$.

Then if $v \in \pi_{n}(\hull_{\cpct}(\{x_{1},x_{2},\dots)) 
\subseteq 
M_{n}(V)$ is given by
\[
  v = p_{n} \sum_{i=1}^{k}\alpha_{k}x_{k}\beta_{k} p_{n}
\]
we consider the map $\phi: \cpct \cBmintensor \cpct \to M_{n}$ 
given by
\[
   w \mapsto \begin{bmatrix} p_{n}\alpha_{1} & \cdots & 
   p_{n}\alpha_{k} & 0 & \cdots \end{bmatrix} w \begin{bmatrix} 
   \beta_{1}p_{n} \\ \vdots \\
   \beta_{k}p_{n}\\ 0 \\ \vdots \end{bmatrix}
\]
so that $v = \phi(\diag(x_{1},x_{2},\ldots)) = 
\phi(\lambda^{-1}(x))$.  But $\phi\comp\lambda^{-1} : 
\cpct \to M_{m}$ is given by
\[
  w \mapsto [a_{i,j}]w[b_{i,j}]
\]
where
\[
  a_{\mu(i,j),\mu(s,t)} = \begin{cases}
    \alpha_{i,j}, & \text{ if $i = 1$, $j \le n$, $s \le k$}\\
    0, & \text{ otherwise,}
    \end{cases}
\]
and
\[
  b_{\mu(i,j),\mu(s,t)} = \begin{cases} \alpha_{i,j}, & \text{ if 
$s 
  = 1$, $t \le n$, $i \le k$}\\
    0, & \text{ otherwise,}
    \end{cases}
\]
and $m = \max \mu(1,\{1,\ldots,n\})$.  So $\phi\comp\lambda^{-1}$ 
is 
an element of $\cb(\cpct,M_{n})$ or, rewriting, an element of 
$M_{n}(\tc)$, and $\|\phi \cmp \lambda^{-1} \|_{\tc} \le 1$.  
So the image of $\hull_{K}(\{x_{1},x_{2},\ldots\})$ sits inside 
$\clos{\ophull(x)}$.  But $\pi_{n}$ is continuous, so if $v_{\nu} 
\in \hull_{K}(\{x_{1},x_{2},\ldots\})$ converges to
\[
  v \in \clos{\hull_{K}(\{x_{1},x_{2},\ldots\})}
\]
then $\pi_{n}(v_{\nu}) \to \pi_{n}(v) \in \clos{\ophull(x)}$.  
Hence $\op{X} \subseteq \clos{\ophull(x)}$.

So we have that across the identification of categories
\[
  \text{$\cpct$-compactness} \implies \text{operator compactness.}
\]

Given that $\cpct$-compact convex sets give rise to matrix convex 
sets, and the fact that the images of unit balls will be convex in the 
appropriate sense, then it is clear that under the equivalence of 
categories $\cpct$-compact mappings become operator compact mappings.

\begin{thm}\label{thm:opapprox2}
  An operator space $V$ has the operator approximation property, if 
  and only if for any operator space $W$ and any operator compact map 
  $\varphi \in \cb(W,V)$, $\varphi$ can be approximated completely 
  uniformly by finite rank maps.
\end{thm}

\begin{proof}
  First assume that $V$ has the operator approximation property.  
  Since $\varphi$ is operator compact, there is an $x \in \cpct(V)$ 
  such that $\op{\varphi}(\op{B}) \subseteq \ophull{x}$, where 
  $\op{B}$ is the matrix unit ball of $W$.  Now we know that for any 
  $\varepsilon > 0$ we can find a finite rank map $\psi \in \cb(V,V)$ 
  such that
  \[
    \|(\id \tensor \psi)(x) - x\| < \varepsilon.
  \] 
  Therefore
  \[
    \|\op{\psi}(v) - v\|_{\infty} < \varepsilon
  \]
  for all $v \in \ophull{x}$ and so
  \[
    \| \psi \cmp \varphi - \varphi \|_{cb} < \varepsilon.
  \]

  The converse follows from Theorem~\ref{thm:C*-approx}.  We observe 
  that if we have a $\cpct$-compact map $\varphi_{\infty}$ then 
  $\varphi$ is an operator compact map and so we can approximate it by 
  finite rank maps $\varphi_{\nu}$ and therefore the finite 
  $\cpct$-rank maps $(\varphi_{\nu})_{\infty}$ approximate 
  $\varphi_{\infty}$.  Theorem~\ref{thm:C*-approx} tells us that 
  $\cpct(V)$ must then satisfy the $\cpct$-approximation property, and 
  Corollary~\ref{cor:approxpropequiv} then gives our result.
\end{proof}

\section{The Strong Operator Approximation Property}
\label{sec:strongopapprox}

Effros and Ruan in their original paper~\cite{effruan:approx} defined 
an operator space $V$ as having the \emph{strong operator 
approximation property} if the identity map $\id: V \to V$ can be 
approximated by finite rank mappings $\varphi_{\nu}$ in the 
\emph{strongly stable point-norm topology}, i.e.  if $\varphi_{\nu} 
\tensor \id \to \id \tensor \id$ point-norm in $V \copmintensor 
\bdd(H)$ for any Hilbert space $H$.  More generally for 
$\varphi_{\nu}$, $\varphi \in \cb(V,W)$ for some $V$, $W$ operator 
spaces, we say that $\varphi_{\nu} \to \varphi$ in the strongly stable 
point-norm topology if $\varphi_{\nu} \tensor \id \to \varphi \tensor 
\id$ in the point-norm topology on $\cb(V \copmintensor \bdd(H),W 
\copmintensor \bdd(H))$ for all Hilbert spaces $H$.  Since $\cpct 
\copmintensor V \hookrightarrow \bdd(\ell^{2}) \copmintensor V$, it is 
easy to see that the strong operator approximation property implies 
the operator approximation property.  Effros and Ruan posed the 
question as to whether or not the strong operator approximation 
property was equivalent to the operator approximation property, 
conjecturing that it was not.

Kirchberg~\cite{kirchberg:opapprox,kirchberg:extensionsexact} showed 
that this is in fact the case.  A sketch of the argument is as 
follows: the strong operator approximation property is the same as the 
general slice map property which, for C*-algebras, implies exactness.  
Extensions of C*-algebras with the operator approximation property 
have the operator approximation property, but there is an extension of 
$\cone C^{*}_{r}(SL(2,\integers))$ by $\cpct$ (both of which have the 
operator approximation property) which is not exact, and so cannot 
have the strong operator approximation property.  He also showed that 
if the operator space is locally reflexive, then they do agree.

So it seems that there should be another, different, notion of 
compactness which corresponds to the strong operator approximation 
property in the same way that operator compactness corresponds to the 
operator approximation property.  An initial guess might be that this 
should be the matrix compactness used in the duality results mentioned 
at the start of Section 2, and this turns out to be correct.  However 
to utilize this we need to reformulate the definition in terms of what 
is essentially an operator space version of total boundedness.

Defining an analogue of total boundedness for operator spaces runs 
into immediate difficulties, since we only really know what balls 
centered at the origin look like.  However we can avoid this by 
noting that a set $K$ is totally bounded if $K$ is bounded and if 
for every $\varepsilon > 0$ there is a finite dimensional subspace 
$V_{\varepsilon}$ such that every point of $K$ lies within 
$\varepsilon$ of a point in $V_{\varepsilon}$.

This implies total boundedness: given any $\varepsilon > 0$, we can 
find a finite dimensional subspace $V_{\varepsilon/3}$ so that for 
every $x \in K$ there is a point $v \in V_{\varepsilon/3}$ such that 
$\| x - v \| < \varepsilon/3$.  Since $V_{\varepsilon/3}$ is finite 
dimensional and $K$ is closed and bounded, we can cover
\[
  S = \{ v \in V_{\varepsilon/3} : d(K,v) < \varepsilon/3 \}^{-}
\]
by finitely many $\varepsilon/3$ balls, centered at $v_{1},\dots,v_{k} 
\in S$.  But then for any $x \in K$, we can find a $v \in S$ so that 
$\| x - v \| \le \varepsilon/3$ and there is an $i$ so that $v_{i}$ 
lies within $\varepsilon/3$ of $v$, and hence $x$ lies within 
$\varepsilon$ of one of the $v_{i}$, and so $K$ is totally bounded.

Conversely, if $K$ is totally bounded, for any $\varepsilon > 0$, 
choose $v_{1},\dots,v_{n}$ be the centers of $\varepsilon$-balls 
which cover $K$. Then
\[
  V_{\varepsilon} = \spn\{v_{1},\dots,v_{n}\}
\]
is a finite dimensional subspace which meets our criterion.

Hans Saar, a student of Wittstock, in his thesis~\cite{saar:thesis} 
implicitly noticed this. He worked with compact maps, but if you 
look at his conditions on the maps, they imply the following about 
the images of unit balls.

\begin{defn}
  A \emph{matrix point} $v = (v_{n})$ in a vector space $V$ is a 
  sequence of points $v_{i} \in M_{n}(V)$ for $i \in \naturals$.

  A matrix set $\op{K}$ in an operator space $V$ is said to be 
  \emph{completely compact} if $\op{K}$ is closed, completely bounded 
  and if for all $\varepsilon > 0$, there exists a finite dimensional 
  subspace $V_{\varepsilon}$ of $V$ such that for every matrix point 
  $(x_{n}) \in \op{K}$ we have a matrix point $(v_{n}) \in 
  \op{V_{\varepsilon}}$, such that $\| x_{n} - v_{n} \|_{n} < 
  \varepsilon$ for all $n$.
\end{defn}

I would like to thank Zhong-Jin Ruan for bringing Saar's work to my 
attention and for providing this definition.

A matrix set which is operator compact is automatically completely 
compact.  First we note that $\clos{\ophull(x)}$ is strongly operator 
compact, since for any $\varepsilon > 0$, we choose $n$ sufficiently 
large that $\| x - p_{n}xp_{n} \| \le \varepsilon$ (we can do this 
since $\cpct(V)$ is the completion of $M_{\infty}(V)$) and let 
$V_{\varepsilon}$ be the subspace spanned by the entries of 
$p_{n}xp_{n}$.  Then given any matrix point $(v_{1},v_{2},\ldots)$ in 
$\ophull(x)$ we let $v_{i} = (\sigma_{i} \tensor \id)(x)$ and so if we 
let $v_{i}' = (\sigma_{i} \tensor \id)(p_{n}xp_{n})$, we have
\[
  \| v_{i}' - v_{i} \| = \|(\sigma_{i} \tensor \id)(p_{n}xp_{n} - 
  x)\| \le \|(p_{n}xp_{n} - x)\| \le \varepsilon.
\]
Taking closures we get that any point in $\clos{\ophull(x)}$ must lie 
within $\varepsilon$ of $V_{\varepsilon}$.  We extend this to 
arbitrary operator compact sets $\op{K} \subseteq \clos{\ophull{x}}$ 
by using the $V_{\varepsilon}$ that you use for $\clos{\ophull{x}}$.  
We will shortly show that the converse is false in general.  We will 
investigate conditions under which the two definitions agree in 
Section 6.

Completely compact sets are also matrix compact, since for each level 
$K_{n}$ of a completely compact set $\op{K}$, and for every 
$\varepsilon > 0$, we have that every element of $K_{n}$ lies within 
$\varepsilon$ of the finite dimensional subspace 
$M_{n}(V_{\varepsilon})$, and so $K_{n}$ is closed and totally 
bounded.  The converse is not true in general.  If we take $V = 
\ell^{2}(\naturals)_{c}$, with standard basis $\{e_{k}\}$ and let 
$X_{n}$ be the unit ball of $M_{n}(\spn\{e_{1},\dots,e_{n}\})$.  Then 
this set is matrix compact, but is not completely compact, since given 
any finite dimensional subspace of $V$, we can find an $e_{m}$ such 
that $d(e_{m},V) > 1 - \varepsilon$.

However if the set $\op{K}$ in question is matrix convex, as will 
suffice for our discussion, then matrix compactness implies complete 
compactness.  To see this we choose a finite $\varepsilon/2$-net 
$\{v_{i}\}$ for $K_{1}$ and let $V_{\varepsilon} = \spn\{v_{i}\}$. Now 
consider the $n$th level of $\op{K}$ and assume that there is an $x 
\in K_{n}$ such that $d(x,M_{n}(V_{\varepsilon})) > \varepsilon$. 
Then in particular, $x$ is at least $\varepsilon$ distant from any 
matrix of the form $v = [v_{i,j}]$ where the $v_{i,j}$ are taken from 
the $\varepsilon/2$-net. But then we have that

We would like first to show that completely compact matrix sets 
are indeed related to the strong operator approximation property.  
To do this we need to introduce the topology of completely uniform 
convergence on completely compact sets:

\begin{defn}
  If $V$ and $W$ are operator spaces, we say that a sequence of maps 
  $\varphi_{\nu} \in \cb(V,W)$ converges to $\varphi \in \cb(V,W)$ 
  \emph{completely uniformly on completely compact sets} if for all 
  completely compact sets $\op{X} \subset V$ and $\varepsilon > 0$ 
  there is an $N$ such that \[\sup \{\| (\varphi_{\nu})_{n}(x) - 
  \varphi_{n}(x) \|_{n} : x \in X_{n},\ n \in \naturals \} < 
  \varepsilon,\ \forall \nu \ge N\]
\end{defn}

Clearly if $\varphi_{\nu} \to \varphi$ completely uniformly on 
completely compact sets, then it converges completely uniformly on 
operator compact sets, and hence in the stable point-norm topology.  I 
am once again indebted to Zhong-Jin Ruan for pointing out the 
following:

\begin{lemma}
  If $\varphi_{\nu}$, $\varphi \in \cb(V,W)$ for some $V$, $W$ 
  operator spaces, and $\varphi_{\nu} \to \varphi$ completely 
  uniformly on completely compact sets then $\varphi_{\nu} \to 
  \varphi$ in the strongly stable point-norm topology.
\end{lemma}

\begin{proof}
  Fix a Hilbert space $H$, some $a \in V \copmintensor \bdd(H)$, and 
  some $\varepsilon > 0$.  Let $a_{\eta} = \sum_{i=1}^{k_{\eta}} 
  v^{\eta}_{i} \tensor \alpha^{\eta}_{i}$ such that $a_{\eta} \to a$.  
  We let $V^{\eta} = \spn \{ v^{\eta}_{1}, \dots, 
  v^{\eta}_{k_{\eta}} \}$, and $P_{n}(\alpha) = p_{n} \alpha p_{n}$, 
  so if $\| a - a_{\eta} \| < \varepsilon$, then
  \[
    \| \id \tensor P_{n} (a) - \id \tensor P_{n} (a_{\eta}) \| < 
    \varepsilon
  \]
  So the matrix set $(\{a_{n} = \id \tensor P_{n} (a)\})$ is completely 
  compact, and hence for all $\nu$ bigger than some $\nu_{0}$, we have 
  that
  \[
	\| \varphi_{\nu} \tensor \id (a) - \varphi \tensor \id (a) \| = 
	\sup_{n} \|(\varphi_{\nu})_{n}(a_{n}) - \varphi_{n}(a_{n}) \| < 
	\delta.
  \]
\end{proof}

We would like to prove a converse result.  What we will prove is 
actually slightly stronger.  We will say that $\varphi_{\nu} \in 
\cb(V,W)$ converges to $\varphi \in \cb(V,W)$ \emph{completely 
uniformly on matrix compact sets} if for all matrix 
compact sets $\op{X} \subset V$ and $\varepsilon > 0$ there is an 
$N$ 
such that
\[
  \sup \{\| (\varphi_{\nu})_{n}(x) - \varphi_{n}(x) \|_{n} : x \in 
  X_{n},\ n \in \naturals \} < \varepsilon,\ \forall \nu \ge N
\]
Since completely compact sets are matrix compact, uniform 
convergence on matrix compact sets implies uniform convergence 
on completely compact sets.

\begin{lemma}
  If $\varphi_{\nu}$, $\varphi \in \cb(V,W)$ for some $V$, $W$ 
  operator spaces, and $\varphi_{\nu} \to \varphi$ in the strongly 
  stable point-norm topology then $\varphi_{\nu} \to \varphi$ 
  completely uniformly on matrix compact sets.
\end{lemma}

\begin{proof}
  Let $\op{X}$ be a matrix compact set in $V$.  Then for each 
  level $X_{n}$, we can find a sequence of points in $x_{n,i} \in 
  M_{n}(V)$ which converge to zero, such that
  \[
    X_{n} \subseteq \clos{\hull(\{x_{n,1},\ldots,x_{n,m},\ldots\})}.
  \]
  Let $x_{n}$ be the matrix in $\bdd(\ell^{2}) \copmintensor M_{n}(V) 
  \isom \bdd(\ell^{2}) \copmintensor V$ given by
  \[
    x_{n} = \diag(x_{n,1},\ldots,x_{n,m},\ldots)
  \]
  and let $x$ be the element of $\bdd(\ell^{2}) \copmintensor 
  \bdd(\ell^{2}) \copmintensor V$ given by
  \[
    x = \diag(x_{1},\dots,x_{n},\dots).
  \]
  Now we observe that any element $v \in 
  \hull(\{x_{n,1},\ldots,x_{n,m},\ldots\})$ can be written as 
  $(\sigma \tensor \id)(x)$ for some $\sigma \in 
  \cb(\bdd(\ell^{2}) \copmintensor \bdd(\ell^{2}),M_{n})$, with 
  $\|\sigma\|_{\cb} \le 1$.  Now since $\varphi_{\nu} \to \varphi$ 
  in the strongly stable point-norm topology, we have that this 
  implies by~\cite{effruan:approx} that
  \[
    \id \tensor \varphi_{\nu} \to \id \tensor \varphi
  \]
  point-norm in $\cb(Z \copmintensor V,Z \copmintensor W)$ for any 
  operator space $Z$, so in particular, we can find an $N$ such that
  \[
	\|(\id \tensor \varphi_{\nu})(x) - (\id \tensor 
	\varphi)(x)\| < \varepsilon
  \]
  for all $\nu \ge N$, and so
  \begin{eqnarray*}
	\|(\varphi_{\nu})_{n}(v) - \varphi_{n}(v)\| & = & \|(\id 
	\tensor \varphi_{\nu})(v) - (\id \tensor \varphi)(v)\| \\
	& = & \|(\sigma \tensor \id)((\id \tensor \varphi_{\nu})(x) - 
	(\id \tensor \varphi)(x))\| \\
	& < & \varepsilon
  \end{eqnarray*}
  for all $\nu \ge N$, and for any $n \in \naturals$ and $v \in 
  \hull(\{x_{n,1},\ldots,x_{n,m},\ldots\})$.  So $\varphi_{\nu}$ 
  converges to $\varphi$ completely uniformly on the matrix set
  \[
    (\hull(\{x_{n,1},\ldots,x_{n,m},\ldots\})),
  \]
  and an $\varepsilon/3$ argument gives us convergence on 
  $\op{X}$.
\end{proof}

So we have proved that these three forms of convergence are 
equivalent. This means that any of these three can be 
substituted for the type of convergence in the strong operator 
approximation property. Kirchberg's result then tells us that 
there is an operator space where the strongly stable point-norm 
topology is different from the stable point-norm topology. Hence 
by the above result and Proposition~\ref{prop:convergence} the 
topology of completely uniform convergence on strongly 
operator compact sets does not agree with the topology of 
completely uniform convergence on operator compact sets. This 
implies:

\begin{lemma}\label{lemma:strongdiffers}
  There is a completely compact set $\op{X}$ in some operator space 
  $V$ such that $\op{X}$ is not operator compact.
\end{lemma}

\begin{proof}
  Let $V$ be a space where the operator approximation property 
  holds, but not the strong operator approximation property. 
  Assume that there was not such $\op{X}$ in this $V$. Then 
  the topology of completely uniform convergence on strongly 
  operator compact sets agrees with the topology of 
  completely uniform convergence on operator compact sets, as 
  there is no difference in the classes of matrix sets, and 
  so we have a contradiction by our previous discussion.
\end{proof}

Turning this around we can give a condition for when the 
operator approximation property will imply the strong 
operator approximation property.

\begin{prop}\label{prop:opappimpliesstropapp}
  If an operator space $V$ satisfies the operator approximation 
  property and every completely compact set is operator compact, 
  then $V$ satisfies the strong operator approximation property.
\end{prop}

\begin{cor}
  If a C*-algebra $A$ satisfies the operator approximation 
  property and every completely compact set is operator compact, 
  then $A$ is exact.
\end{cor}

This may not be a complete characterization, however, since it 
is conceivable that the operator and completely compact sets may be 
different, but the topologies of completely uniform convergence 
agree.

We need to start relating completely compact sets with C*-compact 
sets, just as in the previous section. As one might expect, we have 
the following result:

\begin{prop}
  Let $V$, $W$ be operator spaces and $\varphi_{\nu}$, $\varphi \in 
  \cb(V,W)$. Then the following are equivalent:
  \begin{enumerate}
  \item $\varphi_{\nu} \to \varphi$ in the strongly stable point-norm 
  topology on $\cb(V,W)$.

  \item $\varphi_{\nu} \to \varphi$ completely uniformly on 
  completely compact sets of $V$.

  \item $\varphi_{\nu} \to \varphi$ completely uniformly on 
  matrix compact sets of $V$.

  \item $\id \tensor (\id \tensor \varphi_{\nu}) \to \id \tensor (\id 
  \tensor \varphi)$ point-norm in
  \[
	\bdd_{\bdd(H)}(\seqzero \cBmintensor (\bdd(H) \copmintensor 
	V),\seqzero \cBmintensor (\bdd(H) \copmintensor W))
  \]
  for every $H$.

  \item $\id \tensor \varphi_{\nu} \to \id \tensor \varphi$ uniformly 
  on compact sets of $\bdd(H) \copmintensor V$ for every $H$.

  \item $\id \tensor \varphi_{\nu} \to \id \tensor \varphi$ uniformly 
  on $\bdd(H)$-compact sets of $\bdd(H) \copmintensor V$ for every $H$.
  \end{enumerate}
\end{prop}

\begin{proof}
  (i $\iff$ ii $\iff$ iii): Lemma 5.1 and 5.2

  (i $\implies$ iv): we observe that for any $H$
  \[
    \seqzero \cBmintensor (\bdd(H) \copmintensor V) \isom 
    \seqzero \copmintensor (\bdd(H) \copmintensor V) \isom
    (\seqzero \copmintensor \bdd(H)) \copmintensor V
  \]
  as Banach spaces, by Lemma 4.2 and so if $\varphi_{\nu} \to \varphi$ 
  in the strongly stable point-norm topology, then in particular
  \[
	(\id \tensor \id) \tensor \varphi_{\nu} \to (\id \tensor \id) 
	\tensor \varphi
  \]
  point-norm on $(\seqzero \copmintensor \bdd(H)) \copmintensor V)$.

  (iv $\iff$ v): This is just the classical result.

  (v $\implies$ i): Follows since for any $H$ a point in $\bdd(H) 
  \copmintensor V$ is compact.

  (iv $\iff$ vi): Proposition~\ref{prop:convergence2}.
\end{proof}

\begin{cor}\label{cor:strapproxpropequiv}
  An operator space $V$ has the strong operator approximation property 
  if and only if $\bdd(H) \copmintensor V$ has the 
  $\bdd(H)$-approximation property for every Hilbert space $H$.
\end{cor}

We call a map $\varphi \in \cb(W,V)$ \emph{completely compact} (resp.  
\emph{matrix compact}) if the image of the matrix unit ball of 
$\varphi$ is completely compact (resp.  matrix compact).  
Saar~\cite{saar:thesis} showed that the completely compact maps are a 
closed two-sided ideal under composition.  Since classical compact 
maps are a closed two-sided ideal under composition we see that matrix 
compact maps are also a closed two-sided ideal under composition, for 
if $\varphi_{\nu} \to \varphi$, then the $(\varphi_{\nu})_{n}$ are all 
compact, and hence so is $\varphi_{n} = \lim (\varphi_{\nu})_{n}$ for 
all $n$.

Again we want to use Theorem~\ref{thm:C*-approx} to show that the 
strong operator approximation theorem is equivalent to the density of 
the finite rank maps in the completely compact maps. Again we do this 
by looking at the appropriate C*-compact sets.

We have a bijection between operator spaces $V$ and the collection of 
$\bdd(H)$-operator spaces of the form $\bdd(H) \copmintensor V$ where 
$H$ is any Hilbert space.  To see this all one has to notice is that 
the spaces $M_{n} \copmintensor V$ give the $n$th matrix norm on the 
operator space $V$.

Under this correspondence, a matrix set $\op{X}$ in $V$ is matrix 
convex if and only if the sets
\[
  X_{H} = \{ v \in \bdd(H) \copmintensor V: \gamma^{*} v \gamma \in X_{n} \}
\]
where the $\gamma: \complexes^{n} \to H$ are are isometries are 
C*-convex.  Note that $X_{\complexes^{n}} = X_{n}$.  We say that a 
matrix convex set in $V$ is \emph{$\bdd(H)$-compact} if for every 
Hilbert space $H$, $X_{H}$ is $\bdd(H)$-compact.

\begin{lemma}
  Any $\bdd(H)$-compact matrix convex set is completely compact.
\end{lemma}

\begin{proof}
  Let $X = X_{\ell^{2}}$ and let $X \subseteq \hull_{\bdd(\ell^{2})} \{ 
  x_{i} \}$. Given any $\varepsilon > 0$, we can find
  \[
    y_{i} = \sum_{j=1}^{k_{i}} \alpha_{i,j} \tensor v_{i,j}
  \]
  such that
  \[
    \| x_{i} - y_{i} \| < \varepsilon/3
  \]
  and we know that there is some $N$ such that $\|x_{i}\| < 
  \varepsilon/3$ for all $i > N$. Given this we let
  \[
    V_{\varepsilon} = \spn \{ v_{i,j} : i \le N, 1 \le j \le k_{i} \}.
  \]
  Then for any matrix point $(a_{n}) \in \op{X}$ we regard $a_{n} \in 
  \bdd(\ell^{2})$ by putting it in the top left corner, and observe 
  that it lies in $X$.  Hence given any $\varepsilon$, we can find a 
  $\bdd(\ell^{2})$-convex combination of $x_{i}$ which lies within 
  $\varepsilon/3$ of $a_{n}$, say
  \[
    z_{n} = \sum_{i=1}^{k} \alpha_{i} x_{i} \beta_{i}
  \]
  and then we note that
  \[
    z_{n}' = \sum_{i=1}^{k} \alpha_{i} y_{i} \beta_{i}
  \]
  lies within $\varepsilon/3$ of $z_{n}$ and hence
  \[
    z_{n}'' = \sum_{i=1}^{\min(k,N)} \alpha_{i} y_{i} \beta_{i}
  \]
  satisfies
  \[
    \| a_{n} - z_{n}'' \| \le \| a_{n} - z_{n}' \| + \| 
    \sum_{i=\min(k,N)}^{k} \alpha_{i} y_{i} \beta_{i}
    < \| a_{n} - z_{n} \| + \| z_{n} - z_{n}' \| + \varepsilon/3
    < \varepsilon
  \]
  and moreover $z_{n}'' \in \bdd(\ell^{2}) \copmintensor 
  V_{\varepsilon}$. Hence
  \[
    \| a_{n} - p_{n}z_{n}''p_{n} \| < \varepsilon
  \]
  and so $\op{X}$ is completely compact.
\end{proof}

Thus if $\id \tensor \varphi$ is $\bdd(H)$-compact for every $H$, 
then $\varphi$ is completely compact.

\begin{thm}\label{thm:opapprox3}
  An operator space $V$ has the strong operator approximation 
  property, if and only if for any operator space $W$ and any 
  completely compact map $\varphi \in \cb(W,V)$, $\varphi$ can be 
  approximated completely uniformly by finite rank maps.
\end{thm}

\begin{proof}
  First assume that $V$ has the strong operator approximation 
  property.  Since $\varphi$ is operator compact, there exists a 
  completely compact set $\op{X}$ in $V$ such that 
  $\op{\varphi}(\op{B}) \subseteq \op{X}$, where $\op{B}$ is the 
  matrix unit ball of $W$.  Since $V$ has the strong approximation 
  property, for any $\varepsilon > 0$ there is a finite rank map $\psi 
  \in \cb(V,V)$ such that
  \[
    \|\op{\psi}(v) - v\| < \varepsilon
  \]
  for all $v \in \op{X}$ and so
  \[
    \| \psi \cmp \varphi - \varphi \|_{cb} < \varepsilon.
  \]

  The converse follows from Theorem~\ref{thm:C*-approx}.  We observe 
  that if we have that $\id \tensor \varphi$ is a $\bdd(H)$-compact 
  map for every $H$, then $\varphi$ is a completely compact map and so 
  we can approximate it by finite rank maps $\varphi_{\nu}$ and 
  therefore the finite $\bdd(H)$-rank maps $\id \tensor \varphi_{\nu}$ 
  approximate $\id \tensor \varphi$.  Theorem~\ref{thm:C*-approx} 
  tells us that $\bdd(H) \tensor V$ must then satisfy the 
  $\bdd(H)$-approximation property, and 
  Corollary~\ref{cor:strapproxpropequiv} then gives our result.
\end{proof}

\section{Subcoexact Operator Spaces}

Given the results of the previous section we would like to be able 
to say when the operator compact and completely compact matrix sets 
in an operator space agree.  Heuristically, what happens in the 
classical case is that the two types of compactness agree because in 
a finite dimensional space we can always find a finite sequence of 
points whose convex hull is as ``close'' to the unit ball as we 
like.  By taking better and better finite dimensional approximations 
to our compact set (in the classical sense of the definition of 
complete compactness), and covering more and more closely we can 
build a sequence which converges to zero and whose hull contains our 
original set.

More precisely, what we mean by ``close'' in the above is that any 
point in the convex hull of $\{x_{1},x_{2},\dots,x_{n}\}$ lies in 
the $1+\varepsilon$ ball of the space $V$, or equivalently, the 
map
\[
  \varphi : \ell^{1}_{n} \to V
\]
defined by
\[
  \varphi : e_{i} \mapsto x_{i}
\]
has norm at most $1 + \varepsilon$.  Really we should think of 
this as an isomorphism $\psi$ from $\ell^{1}_{n}/\ker \varphi$ to 
$V$, and it satisfies $\|\psi\|\|\psi^{-1}\| \le 1 + 
\varepsilon$.  Recall that the \emph{Banach-Mazur distance} 
between two finite-dimensional Banach spaces $V$ and $W$ is given 
by
\[
  d_{b}(V,W) = \inf \{ \|\varphi\|\|\varphi^{-1}\|: \varphi \in 
  \bdd(V,W)\text{ is an isomorphism}\}.
\]
Our heuristic can then be restated as saying that for any finite 
dimensional Banach space $V$ and any $\varepsilon > 0$ there is an $n$ 
and a $W \subset \ell^{1}_{n}$ such that $d_{b}(V,\ell^{1}_{n}/W) < 1+ 
\varepsilon$.  As we will see, and as we would expect from the work of 
Pisier~\cite{pisier:exactopsp}, the analogous statement for operator 
spaces is not necessarily true.

To start us down this road, we will look at quotients of 
$\tc_{n}$.

\begin{lemma}
  Let $V = \tc_{n}/W$, that is $V$ is a finite quotient of the $n$ 
  by $n$ trace class operators, then the matrix unit ball of $V$ is 
  operator compact.
\end{lemma}

\begin{proof}
  By Example~\ref{eg:tcball}, we know that the unit ball of 
  $\tc_{n}$ lies inside $\ophull(\tau)$.  So if $v$ is in the 
  $m$th level of the matrix unit ball of $V$, then there is some 
  $\sigma$ in the $m$th level of the matrix unit ball of 
  $\tc_{n}$ such that $v = \pi_{m}(\sigma)$, where $\pi$ is the 
  quotient map.  Then
  \[
	v = (\id \tensor \pi)(\sigma \tensor \id)(\sum_{i,j = 1}^{n} 
	e_{i,j} \tensor \tau_{i,j}) = (\sigma \tensor \id)(\sum_{i,j 
	= 1}^{n} e_{i,j} \tensor \pi(\tau_{i,j}))
  \]
  so the unit ball lies in $\ophull{\pi_{n}(\tau)}$.
\end{proof}

\begin{cor}\label{cor:finitematrixball}
  If $V$ is any finite dimensional operator space, then the matrix 
  unit ball of $V$ is operator compact.
\end{cor}

\begin{proof}
  For $V$ is completely isomorphic (but maybe not completely 
  isometric) via $\varphi$ to $\tc_{n}/W$ for any $n$ and some 
  $W$ depending on $n$.  Without loss of generality we may take 
  $\|\varphi\|_{cb} = 1$.  Hence if we let $x = 
  \varphi^{-1}_{n}(\pi_{n}(\tau))$, we have that the matrix unit 
  ball $\op{X}$ of $V$ has image $\op{\varphi}(\op{X})$ inside 
  the hull of $\pi_{n}(\tau)$ and so $\op{X} \subseteq 
  \ophull(x)$.
\end{proof}

\begin{cor}
  If $V$ is any finite dimensional operator space, then any closed, 
  completely bounded matrix set is operator compact.
\end{cor}

We notice that we have absolutely no control over the norm of $x$ 
in the corollary, but would like to try.  We recall that the 
\emph{completely bounded (or Pisier-) Banach-Mazur} distance 
between two finite dimensional operator spaces $V$ and $W$ is
\[
  d_{cb}(V,W) = \inf \{ \|\varphi\|_{cb}\|\varphi^{-1}\|_{cb}: 
  \varphi \in \cb(V,W)\text{ is a complete isomorphism}\}.
\]
So clearly, no matter how cleverly we were to choose our 
$\varphi$, we must have $\|x\| \ge d_{cb}(V,\tc_{n}/W)$.  But we 
still can vary $n$, so one might expect that by taking larger $n$ 
we might be able to find a better $x$.  We will say that $V$ is 
\emph{$\lambda$-coexact} if for every $\varepsilon > 0$, we can 
find $n$ and a subspace $W \subset \tc_{n}$ such that
\[
  d_{cb}(V,\tc_{n}/W) < \lambda + \varepsilon.
\]
This is a concept dual to Pisier's \emph{$\lambda$-exactness}: a 
finite dimensional space $V$ is $\lambda$-exact if for every 
$\varepsilon > 0$, we can find an $n$ and a subspace $W$ of $M_{n}$ 
such that
\[
  d_{cb}(V,W) < \lambda + \varepsilon.
\]
We note that for finite dimensional $V$, it is immediate from 
duality that $V$ being $\lambda$-coexact implies that $V^{*}$ is 
$\lambda$-exact, and $V$ being $\lambda$-exact implies that 
$V^{*}$ is $\lambda$-coexact.  Since for any given $\lambda$ there 
are operator spaces that are not $\lambda$-exact, there must be 
operator spaces that are not $\lambda$-coexact.

\begin{example}
  Pisier showed that for $r \ge 3$,
  \[
	d_{cb}(\max \ell^{1}_{r},W) \ge \frac{r}{2\sqrt{r - 1}}
  \]
  and so $\max \ell^{1}_{r}$ and $T_{r}$ are not $\lambda$-exact 
  for any $\lambda \le \frac{r}{2\sqrt{r - 1}}$.  Hence $(\max 
  \ell^{1}_{r})^{*} = \min \ell^{\infty}_{r}$ and $M_{r}$ are not 
  $\lambda$-coexact for any $\lambda \le \frac{r}{2\sqrt{r - 1}}$
\end{example}

For infinite dimensional spaces, we say $V$ is $\lambda$-exact if 
every finite-dimensional space is $\lambda$-exact.  We define 
$d_{ex}(V)$ to be the infimum of the $\lambda$ for which $V$ is 
$\lambda$-exact, or $\infty$ if there are no such $\lambda$.  We 
might be tempted to say that $V$ is $\lambda$-coexact if every 
finite-dimensional space is $\lambda$-coexact, but this does not 
quite measure what we are interested in.  Our original aim was to 
control the size of the $x$ for which the unit ball of a finite 
dimensional subspace $W$ of $V$ is contained in 
$\clos{\ophull{x}}$, but clearly now that we are in a larger 
space, so we might be able to get a better $x$ by choosing it 
sitting in $V$, rather than $W$.  However, $x$ must still be in 
the image of a finite dimensional space, so what we want is the 
following:

\begin{defn}
  Let $V$ be an operator space.  We say that $V$ is 
  \emph{$\lambda$-subcoexact} if for every finite dimensional 
  subspace $W$ of $V$ there is another finite dimensional subspace 
  $X$ containing $W$ and which is $\lambda$-coexact.  We define 
  $d_{sce}(V)$ to be the infimum of the $\lambda$ for which $V$ is 
  $\lambda$-subcoexact, or $\infty$ if there is no such $\lambda$.
\end{defn}

This definition give us sufficient control that we can prove the 
following theorem, which is essentially an adaptation of the classical 
result to the operator space situation.  Before we proceed, we need to 
observe that if $V$ is $\lambda$-subcoexact, then any quotient of it 
by a finite dimensional space is also $\lambda$-subcoexact, for if $W 
\subseteq V/X$, then $\pi^{-1}(W)$ is still finite dimensional, and so 
it is contained in a finite-dimensional $\lambda$-coexact subspace $Z$ 
which contains $X$, and clearly $W \subseteq Z/X \isom (\tc/Y)/X$, 
whence the result.

\begin{thm}
  Let $V$ be a $\lambda$-subcoexact operator space for some 
$\lambda < 
  \infty$.  If $\op{K}$ is a completely compact subset of $V$, then 
  $\op{K}$ is operator compact.
\end{thm}

\begin{proof}
  We will construct a sequence of points $x_{i} \in M_{n_{i}}(V)$.  
  Let $\op{K_{1}} = \op{K}$.  Given $\op{K_{i}} \subseteq V/W_{i}$, 
  where $W_{i}$ is finite dimensional, and $\op{K_{i}}$ is completely 
  compact, we choose a finite dimensional subspace $V_{i}$ so that 
  every matrix point $\op{x} \in 2\op{K_{i}}$ has a matrix point 
  $\op{v} \in \op{V_{i}}$ within $4^{-i}$ of it.  Now since $V/W_{i}$ 
  is $\lambda$-coexact, we can find a subspace $U_{i}$ containing 
  $V_{i}$ which is $\lambda$-coexact.  Moreover, there is an $x'_{i}$ 
  sitting inside $M_{n_{i}}(U_{i})$ so that the completely bounded 
  matrix set
  \[
    \{ v \in M_{k}(V_{i}) : d(v,2\op{K_{i}}) < 4^{-i}\} \subseteq 
    \clos{\ophull{x'_{i}}}
  \]
  since this set sits in a finite dimensional space, and since $U_{i}$ 
  is $\lambda$-coexact we may choose this $x'_{i}$ so that
  \[
    \|x'_{i}\| \le \lambda \sup \{\|x\|:x \in 2\op{K_{i}} + 4^{-i} 
+ 
    \varepsilon_{i}\}.
  \]
  We now let $x_{i}$ be an element of $M_{n_{i}}(V)$ in the preimage 
  of $x'_{2}$ such that
  \[
    \|x_{i}\| - \|x'_{i}\| < \varepsilon_{i}.
  \]
  Finally $K_{i+1} = K_{i} + V_{i} \subseteq (V/W_{i})/V_{i} = 
  V/W_{i+1}$.
  
  We note that if $K = \sup\{\|x\|: x \in \op{K}\}$, then 
  $\sup\{\|x\|:x \in \op{K_{i}}\} \le 2^{-i}K$, and so
  \[
    \|x_{k}\| \le 2^{-i+1}\lambda K + 4^{-i} + 2\varepsilon_{i}
  \]
  and so if we choose $\varepsilon_{i} \to 0$, we have that
  \[
    x = \diag(x_{1},x_{2},\dots)
  \]
  is in $\cpct \copmintensor V$.
  
  We claim that $\op{K} \subseteq \clos{\ophull{x}}$ and so is 
  operator compact.  Given $v \in \op{K}$, we can find $v_{1} \in 
  \ophull{x_{1}}$ so that
  \[
    \| 2v - v_{1} \| \le 1/4 + \varepsilon_{1}
  \]
  and then since $2v - v_{1} \in \op{K_{2}}$, we can find $v_{2} \in 
  \ophull{x_{2}}$ so that
  \[
    \| 2(2v - v_{1}) - v_{2} \| \le 4^{-2} + 2\varepsilon_{2}.
  \]
  Repeating this construction, we build a sequence of points $v_{i}$ 
  so that
  \[
	\| 2^{n}v - 2^{n-1}v_{1} - \dots - v_{n} \| \le 4^{-n} + 2 
	\varepsilon_{n}
  \]
  and hence
  \[
	\| v - (2^{-1}v_{1} + 2^{-2}v_{2} + \dots + v_{n})\| \le 2^{-3n} + 
	2^{-n+1} \varepsilon_{2}.
  \]
  But $w_{n} = 2^{-1}v_{1} + 2^{-2}v_{2} + \dots + v_{n}$ is in 
  $\ophull{x}$, since $\sum_{i=1}^{n} 2^{-i} < 1$, and $w_{n} \to v$, 
  so $v \in \clos{\ophull{x}}$.
\end{proof}

We have some immediate corollaries:

\begin{cor}
  Let $V$ be a $\lambda$-subcoexact operator space, where $\lambda < 
  \infty$.  Then $V$ satisfies the operator approximation property if 
  and only if $V$ satisfies the strong operator approximation 
  property.
\end{cor}

\begin{cor}
  If $A$ is a $\lambda$-subcoexact C*-algebra with $\lambda < \infty$, 
  and $A$ satisfies the operator approximation property, then $A$ is 
  exact.
\end{cor}

\begin{cor}
  There are operator spaces which are not $\lambda$-subcoexact for 
  any $\lambda < \infty$.
\end{cor}

\begin{example}
  The trace-class operators $\tc$ are 1-subcoexact, since Effros and 
  Ruan showed in~\cite{effruan:grothendieck} that for every finite 
  dimensional subspace $V$ of $\tc$ and every $\varepsilon > 0$, $V$ 
  is contained in a subspace $W$ for which $d_{cb}(W,\tc_{n}) \le 
  \varepsilon$ for some $n$.  They called spaces which satisfied this 
  property \emph{$\mathcal{J}$ spaces}, hence any $\mathcal{J}$ space 
  is 1-subcoexact.
\end{example}

At this point it is not entirely clear that subcoexactness is a 
natural enough condition to warrant further consideration.  Is it a 
concept which has wide application, or is it merely an ad-hoc 
construction which is useful for this particular application?  It is 
the author's belief that subcoexactness will be an important property, 
and to conclude this paper we include a result which indicates another 
potential application.

\begin{prop}
  If $V$ is an operator space, and $V^{**}$ is subcoexact then $V$ 
is 
  locally reflexive.
\end{prop}

\begin{proof}
  Let $W$ be an arbitrary finite dimensional operator space and 
  $\varphi: W \to V^{**}$.  We want to approximate $\varphi$ in the 
  point-weak-* topology by complete contractions $\varphi_{\nu}: W 
  \to V$.
  
  We first assume that $W = \tc_{n}/X$ for some $X$, so $W^{*} = 
  X^{\perp} \subseteq M_{n}$, and we know that $M_{n}(V^{*}) = 
  \tc_{n}(V)^{*}$, so that
  \[
	\begin{array}{ccc}
	  (M_{n} \copmintensor V)^{*} & \isom & \tc_{n} \copmaxtensor 
	  V^{*} \\
	  \downarrow & & \downarrow \\
	  (X^{\perp} \copmintensor V)^{*} & \isom & W \copmaxtensor 
	  V^{*}
    \end{array}
  \]
  and the bottom row is a complete isometry.  Hence
  \[
	\cb(W,V^{**}) \isom (W \copmaxtensor V^{*})^{*} \isom (X^{\perp} 
	\copmintensor V)^{**} \isom \cb(W,V)^{**}
  \]
  and by the matrix bipolar theorem, we know that the unit ball of 
  $\cb(W,V)$ is weakly dense in the unit ball of $\cb(W,V)^{**}$, so 
  that there is a net of complete contractions $\varphi_{\nu} \in 
  \cb(W,V)$ which converges to $\varphi$, which means that for all $x 
  \in \op{W}$ and $\psi \in \op{V^{*}}$, we have
  \[
	\oplangle \psi, \op{\varphi_{\nu}}(x) \oprangle = \oplangle x 
	\tensor \psi, \varphi_{\nu} \oprangle \to \oplangle x \tensor 
	\psi, \varphi \oprangle = \oplangle \psi, \op{\varphi}(x) 
	\oprangle.
  \]
  Hence $\varphi_{\nu} \to \varphi$ point-weak-*.

  So now if $V^{**}$ is 1-coexact, for any $\varepsilon > 0$, we 
  can find a subspace $Z$ of $V^{**}$ containing $\varphi(W)$, 
  and a complete isomorphism $\theta : Z \to \tc_{n}/X$ for some 
  $n$ and for some $X$, such that 
  $\|\theta\|_{cb}\|\theta^{-1}\|_{cb} \le 1 + \varepsilon$.  So 
  if we let $\|\theta^{-1}\| = 1$, we have that $\theta^{-1}$ 
  can be approximated by $\theta_{\nu} : \tc_{n}/X \to V$, but if 
  we let
  \[
	\psi_{\nu} = \theta_{\nu} \cmp \theta \cmp \varphi : W \to 
	V
  \]
  we have that $\psi_{\nu} \to \varphi$ weak-*, and 
  $\|\psi_{\nu}\| \le 1 + \varepsilon$.  We let 
  $\varphi_{\nu,\varepsilon} = \frac{\psi_{\nu}}{1 + 
  \varepsilon}$.
  
  Hence given any $x \in \op{W}$ and $\psi \in \op{V^{*}}$, for 
  any $\varepsilon > 0$ we have
  \begin{multline*}
	\|\oplangle \psi, \op{\varphi_{\nu,\delta}}(x) \oprangle - 
	\oplangle \psi, \op{\varphi}(x) \oprangle\| \\
	\le \frac{1}{1+\delta}(\|\oplangle \psi, \op{\psi_{\nu,\delta}}(x) 
	\oprangle - \oplangle \psi, \op{\varphi}(x) \oprangle\| + 
	\delta\|\oplangle \psi, \op{\varphi}(x) \oprangle\|)
  \end{multline*}
  and so if $\nu$ is large enough and $\delta$ small enough, 
  this is smaller than $\varepsilon$.  Hence $\varphi_{\nu} \to 
\varphi$ 
  point-weak-* as required.
\end{proof}

\appendix

\section{Operator Spaces and $\cpct$-Operator Spaces}

We now give the proof of the following proposition.

\begin{prop}
  The functors $\cpct$ and $\fr$ implement an equivalence of 
  categories between the category of operator spaces with completely 
  bounded maps and the category of $\cpct$-operator spaces with 
  continuous bi-$\cpct$-linear maps.
\end{prop}

\begin{proof}
  We have an isomorphism of operator spaces $\iota_{V} : \fr(\cpct(V)) 
  \to V$ given by $[v_{i,j}] \to v_{1,1}$ (noting that all other 
  entries are 0): this is clearly an isomorphism of vector spaces, and 
  if $v = [v_{k,l}] \in M_{n}(\fr(\cpct(V)))$, then
  \[
	\| \iota_{V,n}(v) \|_{n} = \| [\iota_{V}([v_{k,l}])] \|_{n} = \| 
	[[v_{k,l}]_{1,1}] \|_{n}
  \]
  but
  \[
	\| v \|_{n} = \| [[v_{k,l}]_{1,1}] \| = \sup_{m} \| 
	[[v_{k,l}]_{1,1}] \|_{m} = \| [[v_{k,l}]_{1,1}] \|_{n}
  \]
  and so $\iota_{V}$ is a completely isometric isomorphism.  
  Furthermore it is a natural transformation, since if $\varphi \in 
  \cb(V,W)$, then we have that
  \[
    \iota_{W}(\fr(\cpct(\varphi))([v_{i,j}])) = 
    \varphi(v_{1,1}) = (\fr(\cpct(\varphi))(\iota_{V}([v_{i,j}])
  \]
  So we have a natural equivalence $\fr \comp \cpct \isom \id$.

  Similarly we have an isomorphism of $\cpct$-operator spaces $\tau : 
  \cpct(\fr(\mathcal{V})) \to \mathcal{V}$ given by $(v_{i,j}) \to 
  \sum e_{i,1} v_{i,j} e_{1,j}$, which again is bijective, is 
  bi-$\cpct$-linear, since if $v = (v_{i,j}) \in 
  \cpct(\fr(\mathcal{V}))$, then
  \[
	\tau(\alpha v \beta) = \sum_{i,j,k,l} e_{i,1} \alpha_{i,k} v_{k,l} 
	\beta_{l,j} e_{1,j} = \sum_{k,l} (\alpha e_{k,1}) v_{k,l} (e_{1,l} 
	\beta) = \alpha \tau(v) \beta,
  \]
  and is an isometric isomorphism since
  \begin{multline*}
	\| \tau(v) \| = \| \sum_{i,j} e_{i,1} v_{i,j} e_{1,j} \| = 
	\sup_{n} \| p_{n} (\sum_{i,j} e_{i,1} v_{i,j} e_{1,j}) p_{n} \| \\
	= \sup_{n} \| \sum_{i,j=1}^{n} e_{i,1} v_{i,j} e_{1,j} \| = \| 
	v\|.
  \end{multline*}
  Again, this is a natural transformation, and so we have a natural 
  equivalence $\cpct \comp \fr \isom \id$.

  Hence the two categories are equivalent.
\end{proof}

\bibliographystyle{plain}

\end{document}